\documentclass[11pt, twoside, a4paper]{article}

\usepackage[T1]{fontenc}

\usepackage{latexsym}
\usepackage{amsmath}
\usepackage{amssymb}
\usepackage{amsfonts}
\usepackage{amsthm}
\usepackage{amscd}
\usepackage{epsfig}
\usepackage{picins}
\usepackage{verbatim}
\usepackage{fancybox}
\usepackage{moreverb}

\newcommand{\T}{{\mathbb T}}

\newcommand{\D}{{\partial}}
\newcommand{\Cl}{{\mathcal C}}

\newtheorem{theorem}{Theorem}
\newtheorem{lemma}[theorem]{Lemma}

\newtheorem{remark}[theorem]{Remark}

\begin{document}

\begin{center}
\Large{\textbf{Two-scale numerical simulation of the weakly compressible 1D isentropic Euler equations}} \\
\end{center}
\begin{center}
\textit{Emmanuel Fr\'enod} \\
L\'EMEL \& LMAM, Universit\'e Europ\'eenne de Bretagne, F-56000 Vannes \\ 
\textit{Alexandre Mouton} \\
IRMA, Universit\'e Louis Pasteur, F-67084 Strasbourg \\ 
\textit{Eric Sonnendr\"ucker} \\
IRMA, Universit\'e Louis Pasteur, F-67084 Strasbourg \\
\end{center}

\section{Introduction}

This paper enters a work programme concerning the development of two-scale numerical methods to solve PDEs with oscillatory singular perturbations linked to physical phenomena. Recently, these methods have been tested on simple problems. For instance, in Ailliot, Fr\'enod and Monbet\cite{Modelling_coastal}, such a method is used to manage the tide oscillation for long term drift forecast of objects in the ocean; in Fr\'enod, Salvarani and Sonnendr\"ucker\cite{Long_time_simulation}, such a method is used to simulate a charged particle beam in a periodic focusing channel. \\

\indent The ultimate goal of this work programme is to propose efficient numerical methods to simulate plasmas submitted to strong magnetic field. Of course, simulations of magnetic confinement fusion are part of this ultimate goal. \\
\indent Before going further towards this ultimate goal, several questions concerning the behaviour of the concerned methods in front of non-linearities and non-smooth solutions need to be investigated. The field of weakly compressible 1D isentropic Euler equations offers a relatively confortable framework in order to tackle these questions. \\
\indent Indeed, the mathematical context established in Grenier\cite{Oscillatory_perturbations}, Klainerman and Majda\cite{Singular_limits,Compressible_incompressible}, Fortenbach, Fr\'enod, Munz and Sonnendr\"ucker\cite{Multiple_scale}, Fr\'enod, Raviart and Sonnendr\"ucker\cite{Two-scale_expansion}, Fr\'enod and Sonnendr\"ucker\cite{Long_time_behavior,Finite_Larmor_radius}, Majda\cite{Majda}, M\'etivier and Schochet\cite{Metivier-Schochet}, Munz\cite{Munz}, Schochet\cite{Symetric_systems_1988,Fast_singular,Symetric_systems_1986,Euler_bounded_domain}, offers a nice marked out way. This is the first motivation of the present paper. \\
\indent A second motivation originates from the fact that robust numerical methods such as finite volumes set out in Godunov\cite{Godunov}, LeVeque\cite{Finite_volume} and Roe\cite{Roe} have been developped in order to solve compressible or incompressible Euler equations. Nevertheless, when they are used to simulate the flow in an experiment with conditions inducing a small Mach number, the CPU time cost of these methods is too high for operational applications because of the very small time step required in order to capture high speed travelling waves that appear in this case. \\

\indent The precise aim of this paper is to develop a two-scale numerical method to solve

\begin{equation}\label{Grenier}
\begin{split}
\D_{t} u^{\epsilon} + \frac{1}{2} \, \D_{x} \big( (u^{\epsilon})^{2}\big) + \D_{x} \big( q^{\epsilon}(\rho^{\epsilon}) \big) + \frac{1}{\epsilon} \, \D_{x} \rho^{\epsilon} &= 0 \, , \\
\D_{t} \rho^{\epsilon} + \D_{x}( \rho^{\epsilon} u^{\epsilon} ) + \frac{1}{\epsilon} \, \D_{x} u^{\epsilon} &= 0 \, , \\ 
u_{| t \, = \, 0}^{\epsilon} &= u_{0} \, ,\\
\rho_{| t \, = \, 0}^{\epsilon} &= \rho_{0} \, ,
\end{split}
\end{equation}
where $(x,t) \in \T^{1} \times [0,T)$. In this model, $\epsilon$ is the Mach number, $u^{\epsilon} = u^{\epsilon}(x,t)$ is the dimensionless speed of the fluid, $\rho^{\epsilon} = \rho^{\epsilon}(x,t)$ is linked to the dimensionless density $\tilde{\rho}^{\epsilon} = \tilde{\rho}^{\epsilon}(x,t)$ by the relation
\begin{equation} \label{rho_rhotilde}
\rho^{\epsilon} = \frac{\tilde{\rho}^{\epsilon}-1}{\epsilon} \, ,
\end{equation}
and
\begin{equation} \label{def_q}
q^{\epsilon}(\rho^{\epsilon}) = \frac{\gamma-2}{2} \, (\rho^{\epsilon})^{2} + \epsilon \, q_{1}^{\epsilon}(\rho^{\epsilon}) \, , 
\end{equation}
where $\gamma$ is the adiabatic coefficient of the fluid and $q_{1}^{\epsilon}(\rho^{\epsilon})$ is regular. This form of the isentropic Euler equations is set out in Grenier\cite{Oscillatory_perturbations} (p. 494) and M\'etivier and Schochet\cite{Metivier-Schochet}. Equation (\ref{Grenier}) is obtained from the classical form of the 1D isentropic Euler equations
\begin{equation} \label{Euler_origin}
\begin{split}
\D_{t}u^{\epsilon} + u^{\epsilon} \, \D_{x}u^{\epsilon} + \frac{1}{\epsilon^{2} \, \tilde{\rho}^{\epsilon}} \, \D_{x}\big(p(\tilde{\rho}^{\epsilon})\big) &= 0 \, ,\\
\D_{t}\tilde{\rho}^{\epsilon} + \D_{x}\big(\tilde{\rho}^{\epsilon} u^{\epsilon}\big) &= 0 \, , \\
u^{\epsilon}(x,0) &= u_{0}(x) \, ,\\
\tilde{\rho}^{\epsilon}(x,0) &= 1 + \epsilon \, \rho_{0}(x) \, ,
\end{split}
\end{equation}
where usually, pressure function is given by $p(\tilde{\rho}^{\epsilon}) = \frac{(\tilde{\rho}^{\epsilon})^{\gamma}}{\gamma}$. Using the relation (\ref{rho_rhotilde}) between $\rho^{\epsilon}$ and $\tilde{\rho}^{\epsilon}$, we obtain
\begin{equation}
\frac{1}{\epsilon^{2} \, \tilde{\rho}^{\epsilon}} \, \D_{x}\big(p(\tilde{\rho}^{\epsilon})\big) = \frac{1}{\epsilon^{2} \, (1+\epsilon \, \rho^{\epsilon})} \, \D_{x} \Big( \frac{(1+\epsilon\,\rho^{\epsilon})^{\gamma}}{\gamma}\Big) = \frac{(1+\epsilon \, \rho^{\epsilon})^{\gamma-2}}{\epsilon} \, \D_{x}\rho^{\epsilon} \, .
\end{equation}
Making and expansion of $(1+\epsilon\,\rho^{\epsilon})^{\gamma-2}$, we obtain
\begin{equation}
\frac{1}{\epsilon^{2} \, \tilde{\rho}^{\epsilon}} \, \D_{x}\big(p(\tilde{\rho}^{\epsilon})\big) = \frac{1}{\epsilon} \, (\D_{x}\rho^{\epsilon})\,\big(1+(\gamma-2)\, \epsilon \, \rho^{\epsilon} + \epsilon^{2} \, q_{1}^{\epsilon}(\rho^{\epsilon})\big) \, ,
\end{equation}
where $q_{1}^{\epsilon}(\rho^{\epsilon})$ is regular, and, by introducing this result in the first equation of (\ref{Euler_origin}), we deduce the first equation of (\ref{Grenier}). In order to obtain the second equation of (\ref{Grenier}), we introduce the relation (\ref{rho_rhotilde}) in the second equation of (\ref{Euler_origin}). \\

\indent In order to achieve this precise aim, inspired by two-scale convergence theory developped in Nguetseng\cite{General_convergence} and Allaire\cite{Homogenization}, we establish a homogenized model describing the mean behaviour of $(u^{\epsilon},\rho^{\epsilon})$. This homogenized model neither contains nor generates high frequency oscillations but only their mean action. We also establish a way to reconstruct the oscillating solution $(u^{\epsilon},\rho^{\epsilon})$ from obtained mean behaviour. The homogenized model and the reconstruction procedure involve an additional variable which is a copy of the time variable and which allows the transfert of informations between the mean behaviour and the high frequency oscillations.

\section{Results}

\setcounter{equation}{0}

\indent In this section, we claim the mathematical results leading to the construction of our two-scale numerical method. \\
\indent After recalling existence and properties of solutions of (\ref{Grenier}), we give a first theorem giving the homogenized model and convergence properties. Then, we set out a numerical approximation for this homogenized model based on Roe's finite volume method which is the two-scale numerical method that allows to solve (\ref{Grenier}). Finally, we give a theorem giving some properties about the convergence of the considered numerical scheme.\\

\indent It is well known (see Grenier\cite{Oscillatory_perturbations}, Klainerman and Majda\cite{Singular_limits}, Majda\cite{Majda}, Metivier and Schochet\cite{Metivier-Schochet}, or Schochet\cite{Symetric_systems_1988,Symetric_systems_1986,Euler_bounded_domain}) that, with initial data satisfying
\begin{equation}
u_{0}, \, \rho_{0} \in H^{s}(\T^{1}) \quad \textnormal{with} \quad s > \frac{3}{2} \, ,
\end{equation}
there exists an existence time $T > 0$, independent of $\epsilon$, such that the system (\ref{Grenier}) admits a unique solution $(u^{\epsilon},\rho^{\epsilon})$ in $\big(\Cl\big(0,T;H^{s}(\T^{1})\big)\big)^{2} \cap \big(\Cl^{1}\big(0,T;H^{s-1}(\T^{1})\big)\big)^{2}$ for any $\epsilon > 0$. Furthermore, the sequence $(u^{\epsilon},\rho^{\epsilon})_{\epsilon \, > \, 0}$ is bounded in $\big(\Cl\big(0,T;H^{s}(\T^{1})\big)\big)^{2} \cap \big(\Cl^{1}\big(0,T;H^{s-1}(\T^{1})\big)\big)^{2}$ for the usual norm. \\

\indent Having this existence result at hand, we can claim the following theorem giving the homogenized model. \\

\begin{theorem}\label{Convergence_epsilon}
If we fix $\epsilon_{0} > 0$ and $s \geq 2$, there exists a constant $K > 0$ such that, for any $\epsilon \in \, ]0,\epsilon_{0}]$, we have
\begin{equation}\label{estimation_epsilon}
\begin{split}
\Bigg(\int_{0}^{T} \Big\| u^{\epsilon}(\cdot,t) - U\big( \cdot,\frac{t}{\epsilon}, t\big) \Big\|_{L^{2}(\T^{1})}^{2} \, dt \Bigg)^{\frac{1}{2}}&\leq K \epsilon \, , \\
\Bigg(\int_{0}^{T} \Big\| \rho^{\epsilon}(\cdot,t) - R\big( \cdot,\frac{t}{\epsilon}, t\big) \Big\|_{L^{2}(\T^{1})}^{2} \, dt \Bigg)^{\frac{1}{2}} &\leq K \epsilon \, ,
\end{split}
\end{equation}
where $U$ and $R$ are defined by
\begin{equation}\label{def_UR}
\begin{split}
U(x,\tau,t) &= F(x-\tau,t) + B(x+\tau,t) + \frac{\overline{u}}{2\pi} \, ,\\
R(x,\tau,t) &= F(x-\tau,t) - B(x+\tau,t) + \frac{\overline{\rho}}{2\pi} \, ,
\end{split}
\end{equation}
with $\displaystyle \overline{u} = \int_{\T^{1}} u_{0}(x) \, dx$, $\displaystyle \overline{\rho} = \int_{\T^{1}}\rho_{0}(x) \, dx$, and $F,B \in \Cl\big(0,T;H^{s}(\T^{1})\big)$ satisfying 
\begin{equation}\label{FB_mass_0}
\int_{T^{1}} F(x,t) \,dx = \int_{\T^{1}} B(x,t) \, dx = 0 \qquad \forall \, t \in [0,T) \, ,
\end{equation}
and 
\begin{equation}\label{def_FB}
\begin{split}
&\D_{t} F + \D_{x} \Big( \frac{\gamma+1}{4} \, F^{2} + \frac{2\overline{u} + (\gamma-1)\overline{\rho}}{4\pi} \, F \Big) = 0 \, , \\
&\D_{t} B + \D_{x} \Big( \frac{\gamma+1}{4} \, B^{2} + \frac{2\overline{u} - (\gamma-1)\overline{\rho}}{4\pi} \, B \Big) = 0 \, , \\
&F_{|_{t\,=\,0}} = \frac{1}{2} \Big( u_{0} + \rho_{0} - \frac{\overline{u}+\overline{\rho}}{2\pi} \Big) \, , \\
&B_{|_{t\,=\,0}} = \frac{1}{2} \Big( u_{0} - \rho_{0} - \frac{\overline{u}-\overline{\rho}}{2\pi} \Big) \, .
\end{split}
\end{equation}
\end{theorem}

\indent A way to interpret this theorem is that if we can compute some good approximations $F_{h}$ and $B_{h}$ of the solutions $F$ and $B$ of (\ref{def_FB}), we are able to reconstruct some functions $U_{h}$ and $R_{h}$ in the same way as we reconstruct $U$ and $R$ in (\ref{def_UR}). Then the obtained functions $(x,t) \mapsto U_{h}(x,\frac{t}{\epsilon},t)$ and $(x,t) \mapsto R_{h}(x,\frac{t}{\epsilon},t)$ are good approximations of $(x,t) \mapsto u^{\epsilon}(x,t)$ and $(x,t) \mapsto \rho^{\epsilon}(x,t)$ respectively. \\
\indent Based on this idea, we build our two-scale numerical method: firstly, we consider a uniform mesh on $\T^{1} \times [0,T]$ with space step $h$ and time step $k$, and we use the notations $x_{i} = ih$ and $t_{n} = nk$. Then we build $F_{h}$ and $B_{h}$ such that
\begin{equation} \label{Roe_FB_1}
F_{h}(x,t) = F_{i}^{n} \textnormal{ and } B_{h}(x,t) = B_{i}^{n} \quad \forall \, (x,t) \in [x_{i-1/2},x_{i+1/2}[ \times [t_{n},t_{n+1}[ \, ,
\end{equation}
where $F_{i}^{n}$ and $B_{i}^{n}$ are computed with Roe's finite volume approximation of (\ref{def_FB}):
\begin{equation} \label{Roe_FB_2}
\begin{split}
F_{i}^{n+1} &= F_{i}^{n} - \frac{k}{2h} \Bigg( \frac{\gamma+1}{4} \big( (F_{i+1}^{n})^{2} - (F_{i-1}^{n})^{2} \big) + \frac{2\overline{u} + (\gamma-1)\overline{\rho}}{4\pi} (F_{i+1}^{n}-F_{i-1}^{n}) \\
&\qquad \qquad \qquad - \Big| \frac{\gamma+1}{4}(F_{i+1}^{n}+F_{i}^{n}) + \frac{2\overline{u}+(\gamma-1)\overline{\rho}}{4\pi} \Big| (F_{i+1}^{n}-F_{i}^{n}) \\
&\qquad \qquad \qquad + \Big| \frac{\gamma+1}{4}(F_{i}^{n}+F_{i-1}^{n}) + \frac{2\overline{u}+(\gamma-1)\overline{\rho}}{4\pi} \Big| (F_{i}^{n}-F_{i-1}^{n}) \Bigg) \, ,
\end{split}
\end{equation}
\begin{equation} \label{Roe_FB_3}
\begin{split}
B_{i}^{n+1} &= B_{i}^{n} - \frac{k}{2h} \Bigg( \frac{\gamma+1}{4} \big( (B_{i+1}^{n})^{2} - (B_{i-1}^{n})^{2} \big) + \frac{2\overline{u} - (\gamma-1)\overline{\rho}}{4\pi} (B_{i+1}^{n}-B_{i-1}^{n}) \\
&\qquad \qquad \qquad - \Big| \frac{\gamma+1}{4}(B_{i+1}^{n}+B_{i}^{n}) + \frac{2\overline{u}-(\gamma-1)\overline{\rho}}{4\pi} \Big| (B_{i+1}^{n}-B_{i}^{n}) \\
&\qquad \qquad \qquad + \Big| \frac{\gamma+1}{4}(B_{i}^{n}+B_{i-1}^{n}) + \frac{2\overline{u}-(\gamma-1)\overline{\rho}}{4\pi} \Big| (B_{i}^{n}-B_{i-1}^{n}) \Bigg) \, .
\end{split}
\end{equation}

\indent Finally, we compute $U_{h}$ and $R_{h}$ by 
\begin{equation} \label{def_UhRh}
\begin{split}
U_{h}(x,\tau,t) &= F_{h}(x-\tau,t) + B_{h}(x+\tau,t) + \frac{\overline{u}}{2\pi} \, ,\\
R_{h}(x,\tau,t) &= F_{h}(x-\tau,t) - B_{h}(x+\tau,t) + \frac{\overline{\rho}}{2\pi} \, ,
\end{split}
\end{equation}
and we have the following convergence result.

\begin{theorem} \label{FV_CV}
If $s \geq 2$, then the approximations $U_{h}(\cdot,\tau,\cdot)$ and $R_{h}(\cdot,\tau,\cdot)$ converge to $U(\cdot,\tau,\cdot)$ and $R(\cdot,\tau,\cdot)$ in $L^{1}\big([0,T) \times \T^{1}\big)$ for any $\tau \in \T^{1}$. Furthermore, if $s \geq 3$, the local truncation errors of the numerical scheme (\ref{Roe_FB_1})-(\ref{Roe_FB_2})-(\ref{Roe_FB_3}) are first order accurate.
\end{theorem}

\begin{remark}
The proof of the convergence result (\ref{estimation_epsilon}) and of Theorem \ref{FV_CV} do not work if $s \in \, ]\frac{3}{2},2[$. However, in this case, the regularity is enough to get 
\begin{equation} \label{wcv_ueps_rhoeps}
\begin{split}
u^{\epsilon} - U(\cdot,\frac{\cdot}{\epsilon},\cdot) & \rightharpoonup 0 \quad weakly-* \, , \\
\rho^{\epsilon} - R(\cdot,\frac{\cdot}{\epsilon},\cdot) & \rightharpoonup 0 \quad weakly-* \, ,
\end{split}
\end{equation}
in $L^{\infty}\big(0,T;H^{s}(\T^{1})\big)$ when $\epsilon \to 0$, with $U$ and $R$ solution of (\ref{def_UR})-(\ref{def_FB}).
\end{remark}

\section{Construction of the two-scale numerical method}

\setcounter{equation}{0}

\subsection{Homogenized model}

\indent We present here the construction of the model (\ref{def_FB}). We recall that we start from the weakly compressible 1D isentropic Euler equations (\ref{Grenier})-(\ref{def_q}). Introducing the functions $f^{\epsilon}, b^{\epsilon}$ defined by
\begin{equation} \label{def_fepsbeps}
\begin{array}{l}
\displaystyle f^{\epsilon}(x,t) = \frac{1}{2}\Big( u^{\epsilon}\big(x+\frac{t}{\epsilon},t\big) + \rho^{\epsilon}\big(x+\frac{t}{\epsilon},t\big) - \frac{\overline{u}+\overline{\rho}}{2\pi} \Big) \, , \\ \\
\displaystyle b^{\epsilon}(x,t) = \frac{1}{2}\Big( u^{\epsilon}\big(x+\frac{t}{\epsilon},t\big) - \rho^{\epsilon}\big(x+\frac{t}{\epsilon},t\big) - \frac{\overline{u}-\overline{\rho}}{2\pi} \Big) \, ,
\end{array}
\end{equation}
we can rewrite the system (\ref{Grenier}): 
\begin{equation} \label{eq_feps}
\begin{split}
& \D_{t} f^{\epsilon}(x,t) + \frac{\epsilon}{2} \, \D_{x} \Bigg( q_{1}^{\epsilon}\Big( f^{\epsilon}(x,t) - b^{\epsilon}\big( x + \frac{2t}{\epsilon} , t \big) + \frac{\overline{\rho}}{2\pi} \Big) \Bigg) \\
& \quad + \D_{x} \Bigg( \frac{\gamma+1}{4} \big(f^{\epsilon}(x,t)\big)^{2} + \frac{\gamma-3}{4} \Big( b^{\epsilon}\big( x + \frac{2t}{\epsilon} , t \big) \Big)^{2} + \frac{2\overline{u} + (\gamma-1) \overline{\rho}}{4\pi} \, f^{\epsilon}(x,t) \\ 
& \qquad \qquad \qquad \qquad + \frac{(3-\gamma)\overline{\rho}}{4\pi} \, b^{\epsilon}\big(x + \frac{2t}{\epsilon} , t\big) + \frac{3-\gamma}{2} \, f^{\epsilon}(x,t) b^{\epsilon}\big(x+\frac{2t}{\epsilon},t\big) \Bigg) = 0 \, ,
\end{split}
\end{equation}
\begin{equation} \label{eq_beps}
\begin{split}
& \D_{t} b^{\epsilon}(x,t) + \frac{\epsilon}{2} \, \D_{x} \Bigg( q_{1}^{\epsilon}\Big( f^{\epsilon}\big(x-\frac{2t}{\epsilon},t\big) - b^{\epsilon}(x,t) + \frac{\overline{\rho}}{2\pi} \Big) \Bigg) \\
& \quad + \D_{x} \Bigg( \frac{\gamma+1}{4} \big(b^{\epsilon}(x,t)\big)^{2} + \frac{\gamma-3}{4} \Big( f^{\epsilon}\big( x - \frac{2t}{\epsilon} , t \big) \Big)^{2} + \frac{2\overline{u} - (\gamma-1) \overline{\rho}}{4\pi} \, b^{\epsilon}(x,t) \\ 
& \qquad \qquad \qquad \qquad - \frac{(3-\gamma)\overline{\rho}}{4\pi} \, f^{\epsilon}\big(x - \frac{2t}{\epsilon} , t\big) + \frac{3-\gamma}{2} \, f^{\epsilon}\big(x-\frac{2t}{\epsilon},t\big) b^{\epsilon}(x,t) \Bigg) = 0 \, , 
\end{split}
\end{equation}
equipped with 
\begin{equation} \label{feps_CI}
f_{|_{t \, = \, 0}}^{\epsilon} = f_{0} = \frac{1}{2}\Big( u_{0} + \rho_{0} - \frac{\overline{u} + \overline{\rho}}{2\pi} \Big) \, ,
\end{equation}
\begin{equation} \label{beps_CI}
b_{|_{t \, = \, 0}}^{\epsilon} = b_{0} = \frac{1}{2}\Big( u_{0} - \rho_{0} - \frac{\overline{u} - \overline{\rho}}{2\pi} \Big) \, .
\end{equation}

Since the sequences $(f^{\epsilon})_{\epsilon \, > \, 0}$ and $(b^{\epsilon})_{\epsilon \, > \, 0}$ are bounded in $L^{\infty}\big(0,T;H^{s}(\T^{1})\big)$, there exist two functions $F$ and $B$ in $L^{\infty}\big(0,T;H^{s}(\T^{1})\big)$ such that for some subsequences always denoted $(f^{\epsilon})_{\epsilon \, > \, 0}$ and $(b^{\epsilon})_{\epsilon \, > \, 0}$, we have
\begin{equation} \label{wcv_fepsF}
f^{\epsilon} \rightharpoonup F \quad \textnormal{weakly-* in $L^{\infty}\big(0,T;H^{s}(\T^{1})\big)$} \, ,
\end{equation}
\begin{equation} \label{wcv_bepsB}
b^{\epsilon} \rightharpoonup B \quad \textnormal{weakly-* in $L^{\infty}\big(0,T;H^{s}(\T^{1})\big)$} \, ,
\end{equation}
for $\epsilon \to 0$.

\indent A first property of $F$ and $B$ is that their averages are equal to 0: to show this, we integrate (\ref{eq_feps})-(\ref{eq_beps}) with respect to $x$ to obtain
\begin{displaymath}
\D_{t} \Bigg(\int_{\T^{1}} f^{\epsilon}(x,t) \, dx \Bigg) = \D_{t} \Bigg(\int_{\T^{1}} b^{\epsilon}(x,t) \, dx \Bigg) = 0 \, .
\end{displaymath}
Then we deduce that, for all $t \in [0,T)$, we have
\begin{equation}
\int_{\T^{1}} f^{\epsilon} \, dx = \int_{\T^{1}} f_{0} \, dx = 0 \, , \quad \int_{\T^{1}} b^{\epsilon} \, dx = \int_{\T^{1}} b_{0} \, dx = 0 \, .
\end{equation}
These results, combined with the fact that $f^{\epsilon}$ and $b^{\epsilon}$ weakly-* converge to $F$ and $B$ respectively, lead to the results (\ref{FB_mass_0}). \\

\indent Furthermore, usual compactness results (for example Lions\cite{Lions}) yield that the functional space 
\begin{equation} \label{spaceU_Aubin-Lions}
\mathcal{U} = \Big\{ g \in L^{\infty}\big(0,T;H^{s}(\T^{1})\big) \, : \, \D_{t}g \in L^{\infty}\big(0,T;H^{s-1}(\T^{1})\big) \Big\} \, ,
\end{equation}
provided with the usual product norm, is compactly embedded in $L^{\infty}\big(0,T;H^{s-1}(\T^{1})\big)$. As a consequence, since $(f^{\epsilon})_{\epsilon \, > \, 0}$ and $(b^{\epsilon})_{\epsilon \, > \, 0}$ are bounded in $\mathcal{U}$, we have
\begin{equation} \label{scv_fepsF}
f^{\epsilon} \stackrel{\epsilon \, \to \, 0}{\longrightarrow} F \quad \textnormal{strongly in $L^{\infty}\big(0,T;H^{s-1}(\T^{1})\big)$} \, ,
\end{equation}
\begin{equation} \label{scv_bepsB}
b^{\epsilon} \stackrel{\epsilon \, \to \, 0}{\longrightarrow} B \quad \textnormal{strongly in $L^{\infty}\big(0,T;H^{s-1}(\T^{1})\big)$} \, .
\end{equation}
Since $f^{\epsilon}$ and $b^{\epsilon}$ are continuous in $t$, we can remark that $F,B \in \Cl\big(0,T;H^{s-1}(\T^{1})\big)$, so the convergence result (\ref{scv_fepsF})-(\ref{scv_bepsB}) is also true in $\Cl\big(0,T;H^{s-1}(\T^{1})\big)$. Having now this convergence result at hand, we can look for the constraint equations on $F$ and $B$. For that, we multiply the equations (\ref{eq_feps})-(\ref{eq_beps}) by a regular function $\varphi$ with compact support on $\T^{1} \times [0,T)$ and we integrate with respect to $x$ and $t$ on $\T^{1} \times [0,T)$. We obtain 
\begin{equation} \label{weq_feps}
\begin{split}
&- \int_{0}^{T} \int_{\T^{1}} f^{\epsilon}(x,t) \, \D_{t}\varphi(x,t) \, dx \, dt - \int_{\T^{1}} f_{0}(x) \, \varphi(x,0) \, dx \\
& \qquad \qquad - \int_{0}^{T} \int_{\T^{1}} \frac{\gamma+1}{4} \big( f^{\epsilon}(x,t)\big)^{2} \, \D_{x} \varphi(x,t) \, dx \, dt \\ 
& \qquad \qquad - \int_{0}^{T} \int_{\T^{1}} \frac{2\overline{u} + (\gamma-1)\overline{\rho}}{4\pi} \, f^{\epsilon}(x,t) \, \D_{x} \varphi(x,t) \, dx \, dt \\
& \qquad \qquad - \int_{0}^{T} \int_{\T^{1}} \frac{\gamma-3}{4} \Big( b^{\epsilon}\big(x+\frac{2t}{\epsilon},t\big)\Big)^{2} \, \D_{x} \varphi(x,t) \, dx \, dt \\
& \qquad \qquad - \int_{0}^{T} \int_{\T^{1}} \frac{(3-\gamma)\overline{\rho}}{4\pi} \, b^{\epsilon}\big(x+\frac{2t}{\epsilon},t\big) \, \D_{x}\varphi(x,t) \, dx \, dt \\
& \qquad \qquad - \int_{0}^{T} \int_{\T^{1}} \frac{3-\gamma}{2} \, f^{\epsilon}(x,t) \, b^{\epsilon}\big(x+\frac{2t}{\epsilon},t\big) \, \D_{x}\varphi(x,t) \, dx \, dt \\
& \qquad \qquad - \frac{\epsilon}{2} \int_{0}^{T} \int_{\T^{1}} q_{1}^{\epsilon}\Big( f^{\epsilon}(x,t) - b^{\epsilon}\big(x+\frac{2t}{\epsilon},t\big) + \frac{\overline{\rho}}{2\pi} \Big) \, \D_{x}\varphi(x,t) \, dx \, dt = 0 \, ,
\end{split}
\end{equation}
and
\begin{equation} \label{weq_beps}
\begin{split}
&- \int_{0}^{T} \int_{\T^{1}} b^{\epsilon}(x,t) \, \D_{t}\varphi(x,t) \, dx \, dt - \int_{\T^{1}} b_{0}(x) \, \varphi(x,0) \, dx \\
& \qquad \qquad - \int_{0}^{T} \int_{\T^{1}} \frac{\gamma+1}{4} \big( b^{\epsilon}(x,t)\big)^{2} \, \D_{x} \varphi(x,t) \, dx \, dt \\ 
& \qquad \qquad - \int_{0}^{T} \int_{\T^{1}} \frac{2\overline{u} - (\gamma-1)\overline{\rho}}{4\pi} \, b^{\epsilon}(x,t) \, \D_{x} \varphi(x,t) \, dx \, dt \\
& \qquad \qquad - \int_{0}^{T} \int_{\T^{1}} \frac{\gamma-3}{4} \Big( f^{\epsilon}\big(x-\frac{2t}{\epsilon},t\big)\Big)^{2} \, \D_{x} \varphi(x,t) \, dx \, dt \\
& \qquad \qquad + \int_{0}^{T} \int_{\T^{1}} \frac{(3-\gamma)\overline{\rho}}{4\pi} \, f^{\epsilon}\big(x-\frac{2t}{\epsilon},t\big) \, \D_{x}\varphi(x,t) \, dx \, dt \\
& \qquad \qquad - \int_{0}^{T} \int_{\T^{1}} \frac{3-\gamma}{2} \, f^{\epsilon}\big(x-\frac{2t}{\epsilon},t\big) \, b^{\epsilon}(x,t) \, \D_{x}\varphi(x,t) \, dx \, dt \\
& \qquad \qquad - \frac{\epsilon}{2} \int_{0}^{T} \int_{\T^{1}} q_{1}^{\epsilon}\Big( f^{\epsilon}\big(x-\frac{2t}{\epsilon},t\big) - b^{\epsilon}(x,t) + \frac{\overline{\rho}}{2\pi} \Big) \, \D_{x}\varphi(x,t) \, dx \, dt = 0 \, .
\end{split}
\end{equation}
\indent Because of (\ref{scv_fepsF}) and (\ref{scv_bepsB}), passing to the limit in the four first terms of (\ref{weq_feps}) and (\ref{weq_beps}) is straightforward. Because of the factor $\frac{\epsilon}{2}$ in front of the last term of (\ref{weq_feps}) and (\ref{weq_beps}), we deduce that these terms converge to 0. \\
\indent In order to find the limit of $\displaystyle \int_{0}^{T} \int_{\T^{1}} b^{\epsilon}\big(x+\frac{2t}{\epsilon},t\big) \, \D_{x} \varphi(x,t) \, dx \, dt$ we firstly make a change of variables
\begin{displaymath}
\begin{split}
& \int_{0}^{T} \int_{\T^{1}} b^{\epsilon}\big(x+\frac{2t}{\epsilon},t\big) \, \D_{x} \varphi(x,t) \, dx \, dt = \int_{0}^{T} \int_{\T^{1}} b^{\epsilon}\big(x,t) \, \D_{x} \varphi\big(x-\frac{2t}{\epsilon},t\big) \, dx \, dt \, ,
\end{split}
\end{displaymath}
then, we write 
\begin{displaymath}
\begin{split}
& \int_{0}^{T} \int_{\T^{1}} b^{\epsilon}\big(x+\frac{2t}{\epsilon},t\big) \, \D_{x} \varphi(x,t) \, dx \, dt \\
&\qquad = \int_{0}^{T} \int_{\T^{1}} \big(b^{\epsilon}\big(x,t)-B(x,t)\big) \, \D_{x} \varphi\big(x-\frac{2t}{\epsilon},t\big) \, dx \, dt + \int_{0}^{T} \int_{\T^{1}} B\big(x,t) \, \D_{x} \varphi\big(x-\frac{2t}{\epsilon},t\big) \, dx \, dt \, .
\end{split}
\end{displaymath}
Then, using the Cauchy-Schwarz inequality, we show that
\begin{equation}
\lim_{\epsilon \, \to \, 0} \int_{0}^{T} \int_{\T^{1}} \big(b^{\epsilon}\big(x,t)-B(x,t)\big) \, \D_{x} \varphi\big(x-\frac{2t}{\epsilon},t\big) \, dx \, dt = 0 \, .
\end{equation}
If we define $\varphi^{\epsilon}(x,t) = \tilde{\varphi}(x,\frac{t}{\epsilon},t)$ where $\tilde{\varphi}(x,\tau,t) = \varphi(x-2\tau,t)$, the sequence $(\varphi^{\epsilon})_{\epsilon \, > \, 0}$ is bounded in $L^{\infty}\big(0,T;H^{s}(\T^{1}))$. Then, there exists a function $\Phi$ in $L^{\infty}\big(0,T;L^{2}(\T^{1};H^{s}(\T^{1}))\big)$ such that
\begin{equation}
\lim_{\epsilon \, \to \, 0} \int_{0}^{T} \int_{\T^{1}} \psi\big(x,\frac{t}{\epsilon},t\big) \, \varphi^{\epsilon}(x,t) \, dx \, dt = \int_{0}^{T} \int_{\T^{1}} \int_{\T^{1}} \psi(x,\tau,t) \, \Phi(x,\tau,t) \, d\tau \, dx \, dt
\end{equation}
for any regular function $\psi$ defined on $\T^{1} \times \T^{1} \times [0,T)$. As described in Allaire\cite{Homogenization}, $\Phi$ is called the two-scale limit of $\varphi^{\epsilon}$. Furthermore, following Marusic-Paloka and Piatnitski\cite{Marusic} and Allaire\cite{Homogenization}, since $\tilde{\varphi}$ is regular and $2\pi$-periodic in $\tau$, it is an easy game to prove that $\Phi = \tilde{\varphi}$ in $L^{\infty}\big(0,T;L^{2}(\T^{1};H^{s}(\T^{1}))\big)$. As a consequence, we have $\Phi(x,\tau,t) = \varphi(x-2\tau,t)$ and
\begin{equation}
\varphi^{\epsilon} \rightharpoonup \int_{\T^{1}} \varphi(x-2\tau,t) \, d\tau \quad \textnormal{weakly-* in $L^{\infty}\big(0,T;H^{s}(\T^{1})\big)$} \, .
\end{equation}
Using these results, we obtain
\begin{equation}
\begin{split}
&\lim_{\epsilon \, \to \, 0} \int_{0}^{T} \int_{\T^{1}} B(x,t) \D_{x} \varphi\big(x-\frac{2t}{\epsilon},t\big) \, dx \, dt \\
&\qquad \qquad = - \int_{0}^{T} \int_{\T^{1}} \int_{\T^{1}} \big(\D_{x} B(x,t)\big) \, \varphi(x-2\tau,t) \, d\tau \, dx \, dt \\
&\qquad \qquad = - \frac{1}{2} \int_{0}^{T} \int_{\T^{1}} B(x,t) \Bigg( \int_{\T^{1}} \D_{\tau} \big(\varphi(x-2\tau,t) \big) \, d\tau \Bigg) \, dx \, dt = 0 \, ,
\end{split}
\end{equation}
and finally
\begin{equation}
\lim_{\epsilon \, \to \, 0} \int_{0}^{T} \int_{\T^{1}} b^{\epsilon}\big(x+\frac{2t}{\epsilon},t\big) \, \D_{x} \varphi(x,t) \, dx \, dt = 0 \, .
\end{equation}

In the same way, we obtain
\begin{align}
\lim_{\epsilon \, \to \, 0} \int_{0}^{T} \int_{\T^{1}} f^{\epsilon}\big(x-\frac{2t}{\epsilon},t\big) \, \D_{x} \varphi(x,t) \, dx \, dt &= 0 \, , \\
\lim_{\epsilon \, \to \, 0} \int_{0}^{T} \int_{\T^{1}} \Big(b^{\epsilon}\big(x+\frac{2t}{\epsilon},t\big)\Big)^{2} \, \D_{x} \varphi(x,t) \, dx \, dt &= 0 \, , \\
\lim_{\epsilon \, \to \, 0} \int_{0}^{T} \int_{\T^{1}} \Big(f^{\epsilon}\big(x-\frac{2t}{\epsilon},t\big)\Big)^{2} \, \D_{x} \varphi(x,t) \, dx \, dt &= 0 \, .
\end{align}

To find the limit of $\displaystyle \int_{0}^{T} \int_{\T^{1}} f^{\epsilon}(x,t) \, b^{\epsilon}\big(x+\frac{2t}{\epsilon},t\big) \, \D_{x} \varphi(x,t) \, dx \, dt$, we notice that
\begin{displaymath}
\begin{split}
&\int_{0}^{T} \int_{\T^{1}} f^{\epsilon}(x,t) \, b^{\epsilon}\big(x+\frac{2t}{\epsilon},t\big) \, \D_{x} \varphi(x,t) \, dx \, dt \\
&\qquad \qquad = \int_{0}^{T} \int_{\T^{1}} \big( f^{\epsilon}(x,t) - F(x,t) \big) \, B\big(x+\frac{2t}{\epsilon},t\big) \, \D_{x} \varphi(x,t) \, dx \, dt \\
& \qquad \qquad \qquad + \int_{0}^{T} \int_{\T^{1}} F(x,t) \, \Big(b^{\epsilon}\big(x+\frac{2t}{\epsilon},t\big) - B\big(x+\frac{2t}{\epsilon},t\big) \Big) \, \D_{x} \varphi(x,t) \, dx \, dt \\
& \qquad \qquad \qquad + \int_{0}^{T} \int_{\T^{1}} \big(f^{\epsilon}(x,t) - F(x,t)\big) \, \Big(b^{\epsilon}\big(x+\frac{2t}{\epsilon},t\big) - B\big(x+\frac{2t}{\epsilon},t\big) \Big) \, \D_{x} \varphi(x,t) \, dx \, dt \\
& \qquad \qquad \qquad + \int_{0}^{T} \int_{\T^{1}} F(x,t) \, B\big(x+\frac{2t}{\epsilon},t\big) \, \D_{x} \varphi(x,t) \, dx \, dt \, .
\end{split}
\end{displaymath}
With Cauchy-Schwarz and H\"older's inequalities, we easily show that
\begin{align}
\int_{0}^{T} \int_{\T^{1}} \big( f^{\epsilon}(x,t) - F(x,t) \big) \, B\big(x+\frac{2t}{\epsilon},t\big) \, \D_{x} \varphi(x,t) \, dx \, dt &\to 0 \, , \\
\int_{0}^{T} \int_{\T^{1}} F(x,t) \, \Big(b^{\epsilon}\big(x+\frac{2t}{\epsilon},t\big) - B\big(x+\frac{2t}{\epsilon},t\big) \Big) \, \D_{x} \varphi(x,t) \, dx \, dt &\to 0 \, , \\
\int_{0}^{T} \int_{\T^{1}} \big(f^{\epsilon}(x,t) - F(x,t)\big) \, \Big(b^{\epsilon}\big(x+\frac{2t}{\epsilon},t\big) - B\big(x+\frac{2t}{\epsilon},t\big) \Big) \, \D_{x} \varphi(x,t) \, dx \, dt &\to 0 \, ,
\end{align}
when $\epsilon \to 0$. If we define the functions $B^{\epsilon}(x,t) = B\big(x+\frac{2t}{\epsilon},t\big)$, we remark that the sequence $(B^{\epsilon})_{\epsilon \, > \, 0}$ is bounded in $L^{\infty}\big(0,T;H^{s-1}(\T^{1})\big)$ and then admits a two-scale limit $\mathfrak{B}$ in $L^{\infty}\big(0,T;L^{\infty}(\T^{1};H^{s-1}(\T^{1}))\big)$ satisfying
\begin{equation}
\mathfrak{B}(x,\tau,t) = B(x+2\tau,t) \qquad \forall \, (x,\tau,t) \in \T^{1} \times \T^{1} \times [0,T) \, ,
\end{equation}
\begin{equation}
B^{\epsilon} \rightharpoonup \int_{\T^{1}} B(x+2\tau,t) \, d\tau = 0 \quad \textnormal{weakly-* in $L^{\infty}\big(0,T;H^{s}(\T^{1})\big)$} \, .
\end{equation}
This means that, for any function $\psi \in \Cl\big(0,T;H^{s-1}(\T^{1})\big)$, we have 
\begin{equation} \label{TS_Beps}
\lim_{\epsilon \, \to \, 0} \int_{0}^{T} \int_{T^{1}} B^{\epsilon}(x,t) \,\psi(x,t) \, dx \, dt = 0 \, .
\end{equation}
Setting then $\psi = F \, \D_{x} \varphi$ in (\ref{TS_Beps}), we obtain
\begin{equation}
\lim_{\epsilon \, \to \, 0} \int_{0}^{T} \int_{\T^{1}} F(x,t) \, B\big(x+\frac{2t}{\epsilon},t\big) \, \D_{x} \varphi(x,t) \, dx \, dt = 0 \, .
\end{equation}
Hence
\begin{equation}
\lim_{\epsilon \, \to \, 0} \int_{0}^{T} \int_{\T^{1}} f^{\epsilon}(x,t) \, b^{\epsilon}\big(x+\frac{2t}{\epsilon},t\big) \, \D_{x} \varphi(x,t) \, dx \, dt = 0 \, .
\end{equation}

The same method also gives
\begin{equation}
\lim_{\epsilon \, \to \, 0} \int_{0}^{T} \int_{\T^{1}} b^{\epsilon}(x,t) \, f^{\epsilon}\big(x-\frac{2t}{\epsilon},t\big) \, \D_{x} \varphi(x,t) \, dx \, dt = 0 \, .
\end{equation}
 
With these results, we can pass to the limit in (\ref{weq_feps})-(\ref{weq_beps}) and get
\begin{equation} \label{weq_F}
\begin{split}
&-\int_{0}^{T} \int_{\T^{1}} F(x,t) \, \D_{t} \varphi(x,t) \, dx \, dt - \int_{\T^{1}} f_{0}(x) \, \varphi(x,0) \,dx \\ &\qquad - \int_{0}^{T} \int_{\T^{1}} \Big( \frac{\gamma+1}{4} \, \big(F(x,t)\big)^{2} + \frac{2\overline{u}+(\gamma-1)\overline{\rho}}{4\pi} \, F(x,t) \Big) \, \D_{x} \varphi(x,t) \, dx \, dt = 0 \, ,
\end{split}
\end{equation}
and
\begin{equation} \label{weq_B}
\begin{split}
& -\int_{0}^{T} \int_{\T^{1}} B(x,t) \, \D_{t} \varphi(x,t) \, dx \, dt - \int_{\T^{1}} b_{0}(x) \, \varphi(x,0) \,dx \\ &\qquad - \int_{0}^{T} \int_{\T^{1}} \Big( \frac{\gamma+1}{4} \, \big(B(x,t)\big)^{2} + \frac{2\overline{u}-(\gamma-1)\overline{\rho}}{4\pi} \, B(x,t) \Big) \, \D_{x} \varphi(x,t) \, dx \, dt = 0 \, .
\end{split}
\end{equation}
Remembering the definition (\ref{feps_CI})-(\ref{beps_CI}) of $f_{0}$ and $b_{0}$, we recognize here the weak formulation of (\ref{def_FB}).

Under the hypothesis of Theorem \ref{Convergence_epsilon} about the initial data in (\ref{def_FB}), we deduce that the solution $(F,B)$ belongs to $\big(\Cl\big(0,T;H^{s}(\T^{1})\big)\big)^{2}$ and is unique in this space (we can apply the Theorem 3.6.1 of Serre\cite{Serre} for example). Finally, from uniqueness, we have all the convergence results above for the whole sequences $(f^{\epsilon})_{\epsilon \, > \, 0}$, $(b^{\epsilon})_{\epsilon \, > \, 0}$ and not only for some subsequences.

\subsection{Convergence of the homogenized model: proof of Theorem \ref{Convergence_epsilon}}

\indent In this section, we will prove inequalities (\ref{estimation_epsilon}). Before going further, we remark that convergence (\ref{wcv_fepsF})-(\ref{wcv_bepsB}) is already proved and that, thanks to (\ref{def_UR}) and (\ref{def_fepsbeps}), these inequalities are equivalent to the theorem below.

\begin{theorem} \label{estimation_epsilon_FB}
For $\epsilon_{0} > 0$ and $s \geq 2$ fixed, there exists a constant $K > 0$ such that
\begin{equation} \label{estimation_epsilon_F}
\Bigg( \int_{0}^{T} \big\| f^{\epsilon}(\cdot,t) - F(\cdot,t) \big\|_{L^{2}(\T^{1})}^{2} \, dt \Bigg)^{\frac{1}{2}} \leq K \epsilon \, ,
\end{equation}
\begin{equation} \label{estimation_epsilon_B}
\Bigg( \int_{0}^{T} \big\| b^{\epsilon}(\cdot,t) - B(\cdot,t) \big\|_{L^{2}(\T^{1})}^{2} \, dt \Bigg)^{\frac{1}{2}} \leq K \epsilon \, ,
\end{equation}
for any $\epsilon \in \, ]0,\epsilon_{0}]$.
\end{theorem}

\indent \textit{Proof of theorem \ref{estimation_epsilon_FB}: } to prove this theorem, in the same spirit of Fr\'enod, Raviart and Sonnendr\"ucker\cite{Two-scale_expansion}, we introduce the functions $\gamma^{\epsilon}$ and $\delta^{\epsilon}$ defined by
\begin{equation}
\begin{split}
\gamma^{\epsilon}(x,t) &= \frac{1}{\epsilon} \big( f^{\epsilon}(x,t) - F(x,t) \big) - W\big(x,\frac{t}{\epsilon},t\big) \, ,\\
\delta^{\epsilon}(x,t) &= \frac{1}{\epsilon} \big(B(x,t) - b^{\epsilon}(x,t) \big) - V\big(x,\frac{t}{\epsilon},t\big) \, ,
\end{split}
\end{equation}
with $W$ and $V$ defined on $\T^{1} \times \T^{1} \times [0,T)$ by
\begin{equation}
\begin{split}
W(x,\tau,t) &= - \frac{\gamma-3}{4} \int_{0}^{\tau} \D_{x} \Big( \big(B(x+2\theta,t)\big)^{2} - 2 F(x,t) B(x+2\theta,t) \\
&\qquad \qquad \qquad \qquad \qquad \qquad \qquad \qquad \qquad \qquad - \frac{\overline{\rho}}{\pi} \, B(x+2\theta,t) \Big) \, d\theta \, ,\\
V(x,\tau,t) &= \frac{\gamma-3}{4} \int_{0}^{\tau} \D_{x} \Big( \big(F(x-2\theta,t)\big)^{2} - 2 F(x-2\theta,t) B(x,t) \\
&\qquad \qquad \qquad \qquad \qquad \qquad \qquad \qquad \qquad \qquad + \frac{\overline{\rho}}{\pi} \, F(x-2\theta,t) \Big) \, d\theta \, .
\end{split}
\end{equation}

Concerning sequences $(\gamma^{\epsilon})_{\epsilon \, \in \, ]0,\epsilon_{0}]}$ and $(\delta^{\epsilon})_{\epsilon \, \in \, ]0,\epsilon_{0}]}$, we have the property below.

\begin{lemma} \label{estimation_Dtgammadelta}
Under the hypothesis of theorem \ref{estimation_epsilon_FB}, there exists a contant $M > 0$, independent of $\epsilon \in \, ]0,\epsilon_{0}]$ and $t \in [0,T)$, such that
\begin{equation} \label{estimation_Dtgamma_Dtdelta}
\begin{split}
&\D_{t} \Big( \big\|\gamma^{\epsilon}(\cdot,t)\big\|_{L^{2}(\T^{1})}^{2} + \big\|\delta^{\epsilon}(\cdot,t)\big\|_{L^{2}(\T^{1})}^{2} \Big) \\
&\qquad \qquad \qquad \qquad \leq M \Big( \frac{3}{2}+\sqrt{2} \Big) \Big( \big\|\gamma^{\epsilon}(\cdot,t)\big\|_{L^{2}(\T^{1})}^{2} + \big\|\delta^{\epsilon}(\cdot,t)\big\|_{L^{2}(\T^{1})}^{2} \Big) + \sqrt{2} M \, ,
\end{split}
\end{equation}
for any $\epsilon \in \, ]0,\epsilon_{0}]$ and $t \in [0,T)$.
\end{lemma}

Having this inequality at hand, we can apply Gronwall's lemma to find a constant $L > 0$, independent of $t$ and $\epsilon$, such that
\begin{equation}
\big\|\gamma^{\epsilon}(\cdot,t)\big\|_{L^{2}(\T^{1})}^{2} + \big\|\delta^{\epsilon}(\cdot,t)\big\|_{L^{2}(\T^{1})}^{2} \leq L \, , \qquad \forall \, (t,\epsilon) \in [0,T) \times ]0,\epsilon_{0}] \, .
\end{equation}

Finally, we define the constant $K$ by
\begin{equation}
K = \sqrt{T} \Bigg[ \sqrt{L} + \max \Big( \sup_{\substack{t \, \in \, [0,T) \\ \tau \, \in \, \T^{1}} } \big\| W(\cdot, \tau, t) \big\|_{H^{s}(\T^{1})} \, , \, \sup_{\substack{t \, \in \, [0,T) \\ \tau \, \in \, \T^{1}} } \big\| V(\cdot, \tau, t) \big\|_{H^{s}(\T^{1})} \Big) \Bigg] \, ,
\end{equation}
and we obtain the inequalities (\ref{estimation_epsilon_F})-(\ref{estimation_epsilon_B}), giving the theorem. \\

\indent Now, to finish the proof of Theorem \ref{estimation_epsilon_FB}, we have to prove Lemma \ref{estimation_Dtgammadelta}. Firstly, we notice that $\D_{t} \gamma^{\epsilon}$ and $\D_{t} \delta^{\epsilon}$ read
\begin{equation} \label{Dt_gamma}
\begin{split}
\D_{t} \gamma^{\epsilon}(x,t) &= \Big( -\alpha \big( f^{\epsilon}(x,t)+F(x,t)\big) + \zeta B\big(x+\frac{2t}{\epsilon},t\big) - \beta_{+} \Big) \, \D_{x} \gamma^{\epsilon}(x,t) \\
&\qquad + \Big( -\alpha \D_{x} \big( f^{\epsilon}(x,t)+F(x,t)\big) + \zeta \D_{x} B\big(x+\frac{2t}{\epsilon},t\big) \Big) \, \gamma^{\epsilon}(x,t) \\
&\qquad + \frac{\zeta}{2} \Big( B\big(x+\frac{2t}{\epsilon},t\big) + b^{\epsilon}\big(x+\frac{2t}{\epsilon},t\big) - 2f^{\epsilon}(x,t) - \frac{\overline{\rho}}{\pi} \Big) \, \D_{x} \delta^{\epsilon}\big(x+\frac{2t}{\epsilon},t\big) \\
&\qquad + \frac{\zeta}{2} \, \D_{x} \Big( B\big(x+\frac{2t}{\epsilon},t\big) + b^{\epsilon}\big(x+\frac{2t}{\epsilon},t\big) - 2f^{\epsilon}(x,t) \Big) \, \delta^{\epsilon}\big(x+\frac{2t}{\epsilon},t\big) \\
&\qquad + \Gamma^{\epsilon}(x,t) \, ,
\end{split}
\end{equation}
\begin{equation} \label{Dt_delta}
\begin{split}
\D_{t} \delta^{\epsilon}(x,t) &= \Big( -\alpha \big( B(x,t)+b^{\epsilon}(x,t)\big) + \zeta F\big(x-\frac{2t}{\epsilon},t\big) - \beta_{-} \Big) \, \D_{x} \delta^{\epsilon}(x,t) \\
&\qquad + \Big( -\alpha \D_{x} \big( B(x,t)+b^{\epsilon}(x,t)\big) + \zeta \D_{x} F\big(x-\frac{2t}{\epsilon},t\big) \Big) \, \delta^{\epsilon}(x,t) \\
&\qquad + \frac{\zeta}{2} \Big( f^{\epsilon}\big(x-\frac{2t}{\epsilon},t\big)+F(x-\frac{2t}{\epsilon},t\big) - 2b^{\epsilon}(x,t) + \frac{\overline{\rho}}{\pi} \Big) \, \D_{x} \gamma^{\epsilon}\big(x-\frac{2t}{\epsilon},t\big) \\
&\qquad + \frac{\zeta}{2} \, \D_{x} \Big( f^{\epsilon}\big(x-\frac{2t}{\epsilon},t\big)+F(x-\frac{2t}{\epsilon},t\big) - 2b^{\epsilon}(x,t) \Big) \, \gamma^{\epsilon}\big(x-\frac{2t}{\epsilon},t\big) \\
&\qquad + \Delta^{\epsilon}(x,t) \, ,
\end{split}
\end{equation}
with $\alpha, \beta_{+}, \beta_{-}, \zeta, \Gamma^{\epsilon}(x,t)$ and $\Delta^{\epsilon}(x,t)$ defined by
\begin{equation} \label{def_alpha_beta_zeta}
\alpha = \frac{\gamma+1}{4}\, , \quad \beta_{+} = \frac{2\overline{u}+(\gamma-1)\overline{\rho}}{4\pi} \, , \quad \beta_{-} = \frac{2\overline{u}-(\gamma-1)\overline{\rho}}{4\pi} \, , \quad \zeta = \frac{\gamma-3}{2} \, ,
\end{equation}
\begin{equation}
\begin{split}
\Gamma^{\epsilon}(x,t) = &-\alpha \D_{x} \Big( W \big(x,\frac{t}{\epsilon},t\big) \big(f^{\epsilon}(x,t)+F(x,t)\big) \Big) - \beta_{+} \D_{x} W\big(x,\frac{t}{\epsilon},t\big) \\
&\, - \zeta \D_{x} \Big( f^{\epsilon}(x,t) V\big(x+\frac{2t}{\epsilon},\frac{t}{\epsilon},t\big) - B\big(x+\frac{2t}{\epsilon},t\big) W\big(x,\frac{t}{\epsilon},t\big) \Big) \\
&\, + \frac{\zeta}{2} \, \D_{x} \Bigg( V\big(x+\frac{2t}{\epsilon},\frac{t}{\epsilon},t\big) \Big( B\big(x+\frac{2t}{\epsilon},t\big)+b^{\epsilon}\big(x+\frac{2t}{\epsilon},t\big)\Big) \Bigg) \\
&\, - \frac{\zeta\overline{\rho}}{2\pi} \, \D_{x} V\big(x+\frac{2t}{\epsilon},\frac{t}{\epsilon},t\big) - \frac{1}{2} \, \D_{x} \Bigg( q_{1}^{\epsilon} \Big( \rho^{\epsilon}\big(x+\frac{t}{\epsilon},t\big)\Big) \Bigg) - \big(\D_{t} W\big) \big(x,\frac{t}{\epsilon},t\big)
\end{split}
\end{equation}
and 
\begin{equation}
\begin{split}
\Delta^{\epsilon}(x,t) = &-\alpha \D_{x} \Big( V\big(x,\frac{t}{\epsilon},t\big)\big(B(x,t)+b^{\epsilon}(x,t) \Big) - \beta_{-} \D_{x} V\big(x,\frac{t}{\epsilon},t\big) \\
&\, -\zeta \D_{x} \Big( W\big(x-\frac{2t}{\epsilon},\frac{t}{\epsilon},t\big) b^{\epsilon}(x,t) - F\big(x-\frac{2t}{\epsilon},t) V\big(x,\frac{t}{\epsilon},t\big) \Big) \\
&\, + \frac{\zeta}{2} \, \D_{x} \Bigg( W\big(x-\frac{2t}{\epsilon},\frac{t}{\epsilon},t\big) \Big( f^{\epsilon}\big(x-\frac{2t}{\epsilon},t\big) + F\big(x-\frac{2t}{\epsilon},t\big) \Big) \Bigg) \\
&\, + \frac{\zeta\overline{\rho}}{2\pi} \, \D_{x} W\big(x-\frac{2t}{\epsilon},\frac{t}{\epsilon},t\big)  + \frac{1}{2} \, \D_{x} \Bigg( q_{1}^{\epsilon} \Big( \rho^{\epsilon} \big( x-\frac{t}{\epsilon},t\big) \Big) \Bigg) - \big(\D_{t} V\big) \big(x,\frac{t}{\epsilon},t\big) \, .
\end{split}
\end{equation}

Since the sequences $(f^{\epsilon})_{\epsilon \, \in \, ]0,\epsilon_{0}]}$ and $(b^{\epsilon})_{\epsilon \, \in \, ]0,\epsilon_{0}]}$ are bounded in $\Cl\big(0,T;H^{s}(\T^{1})\big)$ with $s \geq 2$, $F$ and $B$ are in $\Cl\big(0,T;H^{s}(\T^{1})\big)$, and $q_{1}^{\epsilon}$ is regular, we can find some constants $C_{i} > 0$ ($i = 1,\dots,8$), independent of $\epsilon$ and $t$, such that
\begin{equation} \label{C1}
\begin{split}
2 \int_{\T^{1}} \Big( -\alpha \big( f^{\epsilon}(x,t)+F(x,t)\big) - \beta_{+} + \zeta B\big(x+\frac{2t}{\epsilon},t\big) \Big) \, \gamma^{\epsilon}&(x,t) \, \D_{x} \gamma^{\epsilon}(x,t) \, dx \\
&\leq C_{1} \big\| \gamma^{\epsilon}(\cdot,t)\big\|_{L^{2}(\T^{1})}^{2} \, ,
\end{split}
\end{equation}
\begin{equation} \label{C2}
\begin{split}
2 \int_{\T^{1}} \Big( -\alpha \big( b^{\epsilon}(x,t)+B(x,t)\big) - \beta_{-} + \zeta F\big(x-\frac{2t}{\epsilon},t\big) \Big) \, \delta^{\epsilon}&(x,t) \, \D_{x} \delta^{\epsilon}(x,t) \, dx \\
&\leq C_{2} \big\| \delta^{\epsilon}(\cdot,t)\big\|_{L^{2}(\T^{1})}^{2} \, ,
\end{split}
\end{equation}
\begin{equation} \label{C3}
\begin{split}
2 \int_{\T^{1}} \Big( -\alpha \D_{x} \big( f^{\epsilon}(x,t)+F(x,t)\big) +\zeta \D_{x} B\big(x+\frac{2t}{\epsilon},t\big) \Big) \,& \big(\gamma^{\epsilon}(x,t)\big)^{2} \,dx \\
&\leq C_{3} \big\| \gamma^{\epsilon}(\cdot,t)\big\|_{L^{2}(\T^{1})}^{2} \, ,
\end{split}
\end{equation}
\begin{equation} \label{C4}
\begin{split}
2 \int_{\T^{1}} \Big( -\alpha \D_{x} \big( b^{\epsilon}(x,t)+B(x,t)\big) + \zeta \D_{x} F\big(x-\frac{2t}{\epsilon},t\big) \Big) \,& \big(\delta^{\epsilon}(x,t)\big)^{2} \,dx \\
&\leq C_{4} \big\| \delta^{\epsilon}(\cdot,t)\big\|_{L^{2}(\T^{1})}^{2} \, ,
\end{split}
\end{equation}
\begin{equation} \label{C5}
\begin{split}
2 \int_{\T^{1}} \frac{\zeta}{2} \, \D_{x} \Big( B\big(x+\frac{2t}{\epsilon},t\big) + b^{\epsilon}\big(x+\frac{2t}{\epsilon},t\big) - 2f^{\epsilon}&(x,t) \Big) \, \delta^{\epsilon}\big(x+\frac{2t}{\epsilon},t\big) \, \gamma^{\epsilon}(x,t) \, dx \\
&\leq C_{5} \big\| \gamma^{\epsilon}(\cdot,t)\big\|_{L^{2}(\T^{1})} \big\| \delta^{\epsilon}(\cdot,t)\big\|_{L^{2}(\T^{1})} \, ,
\end{split}
\end{equation}
\begin{equation} \label{C6}
\begin{split}
2 \int_{\T^{1}} \frac{\zeta}{2} \, \D_{x} \Big( F\big(x-\frac{2t}{\epsilon},t\big) + f^{\epsilon}\big(x-\frac{2t}{\epsilon},t\big) - 2b^{\epsilon}&(x,t) \Big) \, \gamma^{\epsilon}\big(x-\frac{2t}{\epsilon},t\big) \, \delta^{\epsilon}(x,t) \, dx \\
&\leq C_{6} \big\| \gamma^{\epsilon}(\cdot,t)\big\|_{L^{2}(\T^{1})} \big\| \delta^{\epsilon}(\cdot,t)\big\|_{L^{2}(\T^{1})} \, ,
\end{split}
\end{equation}
\begin{equation} \label{C7}
2 \int_{\T^{1}} \Gamma^{\epsilon}(x,t) \, \gamma^{\epsilon}(x,t) \, dx \leq C_{7} \big\| \gamma^{\epsilon}(\cdot,t)\big\|_{L^{2}(\T^{1})} \, ,
\end{equation}
\begin{equation} \label{C8}
2 \int_{\T^{1}} \Delta^{\epsilon}(x,t) \, \delta^{\epsilon}(x,t) \, dx \leq C_{8} \big\| \delta^{\epsilon}(\cdot,t)\big\|_{L^{2}(\T^{1})} \, .
\end{equation}
On the other hand, making in the first integral the change of variable $ x \mapsto x + \tau$ and in the second one $x \mapsto x - \tau$, rearranging the terms and using the definitions of $\gamma^{\epsilon}$ and $\delta^{\epsilon}$, we have
\begin{displaymath}
\begin{split}
& \zeta \int_{\T^{1}} \Big( b^{\epsilon}(x+2\tau,t) + B(x+2\tau,t) - 2 f^{\epsilon}(x,t) - \frac{\overline{\rho}}{\pi} \Big) \, \D_{x} \delta^{\epsilon}(x+2\tau,t) \, \gamma^{\epsilon}(x,t) \, dx \\
&\quad + \zeta \int_{\T^{1}} \Big( f^{\epsilon}(x-2\tau,t) + F(x-2\tau,t) - 2 b^{\epsilon}(x,t) + \frac{\overline{\rho}}{\pi} \Big) \, \D_{x} \gamma^{\epsilon}(x-2\tau,t) \, \delta^{\epsilon}(x,t) \, dx \\
&\qquad = \zeta \int_{\T^{1}} \big( b^{\epsilon}(x+\tau,t) - B(x+\tau,t)\big) \gamma^{\epsilon}(x-\tau,t) \, \D_{x} \delta^{\epsilon}(x+\tau,t) \, dx \\
&\qquad \quad + 2\zeta \int_{\T^{1}} \big( F(x-\tau,t) - f^{\epsilon}(x-\tau,t) \big) \gamma^{\epsilon}(x-\tau,t) \, \D_{x} \delta^{\epsilon}(x+\tau,t) \, dx \\
&\qquad \quad - 2\zeta \int_{\T^{1}} \Big( F(x-\tau,t) - B(x+\tau,t) + \frac{\overline{\rho}}{2\pi} \Big) \gamma^{\epsilon}(x-\tau,t) \, \D_{x} \delta^{\epsilon}(x+\tau,t) \, dx \\
&\qquad \quad + \zeta \int_{\T^{1}} \big( f^{\epsilon}(x-\tau,t) - F(x-\tau,t)\big) \delta^{\epsilon}(x+\tau,t) \, \D_{x} \gamma^{\epsilon}(x-\tau,t) \, dx \\
&\qquad \quad + 2\zeta \int_{\T^{1}} \big( B(x+\tau,t) - b^{\epsilon}(x+\tau,t) \big) \delta^{\epsilon}(x+\tau,t) \, \D_{x} \gamma^{\epsilon}(x-\tau,t) \, dx \\
&\qquad \quad + 2\zeta \int_{\T^{1}} \Big( F(x-\tau,t) - B(x+\tau,t) + \frac{\overline{\rho}}{2\pi} \Big) \delta^{\epsilon}(x+\tau,t) \, \D_{x} \gamma^{\epsilon}(x-\tau,t) \, dx \\
&\qquad = \zeta \int_{\T^{1}} \epsilon \big( \D_{x} \gamma^{\epsilon}(x-\tau,t) - \D_{x} \delta^{\epsilon}(x+\tau,t) \big) \, \gamma^{\epsilon}(x-\tau,t) \, \delta^{\epsilon}(x+\tau,t) \, dx \\
&\qquad \quad + \zeta \int_{\T^{1}} \epsilon \D_{x} \gamma^{\epsilon}(x-\tau,t) \big( W(x-\tau,\tau,t) + 2V(x+\tau,\tau,t) \big) \, \delta^{\epsilon}(x+\tau,t) \, dx \\
&\qquad \quad - \zeta \int_{\T^{1}} \epsilon \D_{x} \delta^{\epsilon}(x+\tau,t) \big( V(x+\tau,\tau,t) + 2W(x-\tau,\tau,t) \big) \, \gamma^{\epsilon}(x-\tau,t) \, dx \\
&\qquad \quad - 2\zeta \int_{\T^{1}} \epsilon \D_{x} \delta^{\epsilon}(x+\tau,t) \big( \gamma^{\epsilon}(x-\tau,t) \big)^{2} \, dx \\
&\qquad \quad + 2\zeta \int_{\T^{1}} \epsilon \D_{x} \gamma^{\epsilon}(x-\tau,t) \big( \delta^{\epsilon}(x+\tau,t) \big)^{2} \, dx \\
&\qquad \quad + 2\zeta \int_{\T^{1}} R(x,\tau,t) \big( \delta^{\epsilon}(x+\tau,t) \D_{x} \gamma^{\epsilon}(x-\tau,t) - \gamma^{\epsilon}(x-\tau,t) \D_{x} \delta^{\epsilon}(x+\tau,t) \big) \, dx \, .
\end{split}
\end{displaymath}
The sequences $(\epsilon \gamma^{\epsilon})_{\epsilon \, \in \, ]0,\epsilon_{0}]}$ and $(\epsilon \delta^{\epsilon})_{\epsilon \, \in \, ]0,\epsilon_{0}]}$ are bounded in $\Cl\big(0,T;H^{s}(\T^{1})\big)$, so we can find some constants $C_{i} > 0$ ($i = 9,\dots,13$), independent of $\epsilon$ and $t$, such that
\begin{equation} \label{C9}
\begin{split}
\zeta \int_{\T^{1}} \epsilon \big( \D_{x} \gamma^{\epsilon}(x-\tau,t) - \D_{x} \delta^{\epsilon}(x+\tau,t) \big) \, &\gamma^{\epsilon}(x-\tau,t) \, \delta^{\epsilon}(x+\tau,t) \, dx \\
&\leq C_{9} \big\| \gamma^{\epsilon}(\cdot,t) \big\|_{L^{2}(\T^{1})} \big\| \gamma^{\epsilon}(\cdot,t) \big\|_{L^{2}(\T^{1})} \, ,
\end{split}
\end{equation}
\begin{equation} \label{C10}
\begin{split}
\zeta \int_{\T^{1}} \epsilon \D_{x} \gamma^{\epsilon}(x-\tau,t) \big( W(x-\tau,\tau,t) + 2V(x+\tau,\tau,t) \big) \, &\delta^{\epsilon}(x+\tau,t) \, dx \\
&\leq C_{10} \big\| \delta^{\epsilon}(\cdot,t)\big\|_{L^{2}(\T^{1})} \, , 
\end{split}
\end{equation}
\begin{equation} \label{C11}
\begin{split}
\zeta \int_{\T^{1}} \epsilon \D_{x} \delta^{\epsilon}(x+\tau,t) \big( V(x+\tau,\tau,t) + 2W(x-\tau,\tau,t) \big) \, &\gamma^{\epsilon}(x-\tau,t) \, dx \\
&\leq C_{11} \big\| \gamma^{\epsilon}(\cdot,t)\big\|_{L^{2}(\T^{1})} \, ,
\end{split}
\end{equation}
\begin{equation} \label{C12}
2\zeta \int_{\T^{1}} \epsilon \D_{x} \delta^{\epsilon}(x+\tau,t) \big( \gamma^{\epsilon}(x-\tau,t) \big)^{2} \, dx \leq C_{12} \big\| \gamma^{\epsilon}(\cdot,t)\big\|_{L^{2}(\T^{1})}^{2} \, ,
\end{equation}
\begin{equation} \label{C13}
2\zeta \int_{\T^{1}} \epsilon \D_{x} \gamma^{\epsilon}(x-\tau,t) \big( \delta^{\epsilon}(x+\tau,t) \big)^{2} \, dx \leq C_{13} \big\| \delta^{\epsilon}(\cdot,t)\big\|_{L^{2}(\T^{1})}^{2} \, .
\end{equation}
For the last integral, we simply compute:
\begin{equation} \label{C14}
\begin{split}
& 2\zeta \int_{\T^{1}} R(x,\tau,t) \big( \delta^{\epsilon}(x+\tau,t) \D_{x} \gamma^{\epsilon}(x-\tau,t) - \gamma^{\epsilon}(x-\tau,t) \D_{x} \delta^{\epsilon}(x+\tau,t) \big) \, dx \\
&\qquad = 2\zeta \int_{\T^{1}} \D_{\tau} \big( R(x,\tau,t) \big) \delta^{\epsilon}(x+\tau,t) \, \gamma^{\epsilon}(x-\tau,t) \, dx \\
&\qquad \quad - 2\zeta \D_{\tau} \Bigg( \int_{\T^{1}} R(x,\tau,t) \, \delta^{\epsilon}(x+\tau,t) \, \gamma^{\epsilon}(x-\tau,t) \, dx \Bigg) \\
&\qquad = 2\zeta \int_{\T^{1}} \D_{\tau} \big( R(x,\tau,t) \big) \delta^{\epsilon}(x+\tau,t) \, \gamma^{\epsilon}(x-\tau,t) \, dx \\
&\qquad = -2\zeta \int_{\T^{1}} \delta^{\epsilon}(x+\tau,t) \, \gamma^{\epsilon}(x-\tau,t) \, \D_{x}U(x,\tau,t) \, dx \\
&\qquad \leq C_{14} \big\| \gamma^{\epsilon}(\cdot,t)\big\|_{L^{2}(\T^{1})} \big\| \delta^{\epsilon}(\cdot,t)\big\|_{L^{2}(\T^{1})}
\end{split}
\end{equation}
for a constant $C_{14} > 0$ independent of $\epsilon$ and $t$.

Having (\ref{C1})-(\ref{C14}) at hand, we multiply (\ref{Dt_gamma}) by $\gamma^{\epsilon}(x,t)$ and (\ref{Dt_delta}) by $\delta^{\epsilon}(x,t)$ to deduce
\begin{equation} \label{last_estimate}
\begin{split}
\D_{t} \Big( \big\| \gamma^{\epsilon}(\cdot,t)\big\|_{L^{2}(\T^{1})}^{2} + \big\| \delta^{\epsilon}(\cdot,t)\big\|_{L^{2}(\T^{1})}^{2} \Big) &\leq M \Big( \big\| \gamma^{\epsilon}(\cdot,t)\big\|_{L^{2}(\T^{1})}^{2} + \big\| \delta^{\epsilon}(\cdot,t)\big\|_{L^{2}(\T^{1})}^{2} \\
&\qquad \qquad + \big\| \gamma^{\epsilon}(\cdot,t)\big\|_{L^{2}(\T^{1})} + \big\| \delta^{\epsilon}(\cdot,t)\big\|_{L^{2}(\T^{1})} \\
&\qquad \qquad + \big\| \gamma^{\epsilon}(\cdot,t)\big\|_{L^{2}(\T^{1})} \big\| \delta^{\epsilon}(\cdot,t)\big\|_{L^{2}(\T^{1})} \Big) \, ,
\end{split}
\end{equation}
with $M = \max( C_{1}+C_{3}+C_{12}, C_{2}+C_{4}+C_{13}, C_{5}+C_{6}+C_{9}+C_{14},C_{7}+C_{11},C_{8}+C_{10})$. Inequality (\ref{estimation_Dtgamma_Dtdelta}) is directly obtained from (\ref{last_estimate}), ending the proof of lemma \ref{estimation_Dtgammadelta} and then the proof of theorem \ref{estimation_epsilon_FB}. $\square$ \\ \\ \\

\subsection{Properties of the finite volume scheme: proof of the Theorem 2}

Thanks to (\ref{def_UR}) and (\ref{def_UhRh}), convergence of $\big(U_{h}(\cdot,\tau,\cdot), R_{h}(\cdot,\tau,t)\big)$ to $\big(U(\cdot,\tau,\cdot), R(\cdot,\tau,t)\big)$ for any $\tau \in \T^{1}$ is equivalent to the convergence of $(F_{h},B_{h})$ to $(F,B)$. \\
\indent Furthermore, the equations satisfied by $F$ and $B$ are of the form
\begin{equation} \label{eq_q}
\D_{t} q + \D_{x} \big( f(q) \big) = 0 \, , \quad q(x,0) = q_{0}(x) \, , \quad \int_{\T^{1}} q(x,t) \, dx = 0 \, , \quad \forall \, t \, ,
\end{equation}
with $f(q) = \alpha q^{2} + \beta q$, $\alpha$ and $\beta = \beta_{\pm}$ defined in (\ref{def_alpha_beta_zeta}). In the same way, the numerical method (\ref{Roe_FB_1})-(\ref{Roe_FB_2})-(\ref{Roe_FB_3}) is of the form
\begin{equation} \label{Roe_Q}
\begin{split}
Q_{h}(x,t) &= Q_{i}^{n} \qquad \forall \, [x_{i-1/2},x_{i+1/2}[ \times [t_{n},t_{n+1}[  \, , \\
Q_{i}^{n+1} &= Q_{i}^{n} - \frac{k}{h} \big( \mathcal{F}(Q_{i+1}^{n},Q_{i}^{n}) - \mathcal{F}(Q_{i}^{n},Q_{i-1}^{n}) \big) \, , \\
\mathcal{F}(Q_{i},Q_{i-1}) &= \frac{1}{2} \big( f(Q_{i}) + f(Q_{i-1}) \big) - \frac{1}{2} \big| \alpha (Q_{i}+Q_{i-1}) + \beta \big| (Q_{i}-Q_{i-1}) \, ,
\end{split}
\end{equation}
with a space step $h = \cfrac{2\pi}{N_{x}+1}$ and a time step $k$. \\

We notice that the CFL condition at the $n$-th time step is given by
\begin{equation} \label{CFL}
\frac{k}{h} \, \max_{0 \, \leq \, i \, \leq \, N_{x}} \big| \alpha(Q_{i}^{n}+Q_{i-1}^{n}) + \beta \big|  = \nu \leq 1 \, .
\end{equation}

With these notations, proving Theorem 2 is equivalent to proving the theorem below.

\begin{theorem}
If $s \geq 2$, the approximation $Q_{h}$ converges to the solution $q$ of (\ref{eq_q}) in $L^{1} \big( [0,T) \times \T^{1} \big)$ norm. Furthermore, if $s \geq 3$, the local truncation error of the numerical scheme (\ref{Roe_Q}) is first order accurate.
\end{theorem}

\indent \textit{Proof of theorem 5:} the \textit{total variation} of a function $q \in L^{1}\big( [0,T) \times \T^{1} \big)$ is defined by
\begin{equation}
TV_{T}(q) = \limsup_{\eta \, \to \, 0} \frac{1}{\eta} \int_{0}^{T} \int_{\T^{1}} \Big[ \big| q(x+\eta,t) - q(x,t) \big| + \big| q(x,t+\eta) - q(x,t) \big| \Big] \, dx \, dt \, .
\end{equation}
For an approximation $Q_{h}$ computed with the numerical method (\ref{Roe_Q}), we have
\begin{equation}
TV_{T}(Q_{h}) = \sum_{n \, = \, 0}^{T/k} \big[ k \, TV(Q^{n}) + \| Q^{n+1}-Q^{n} \|_{1} \big] \, ,
\end{equation}
with $TV(Q)$ and $\|Q\|_{1}$ defined by
\begin{equation}
TV(Q) = \sum_{i\,=\,0}^{N_{x}} |Q_{i}-Q_{i-1}| \quad \textnormal{and} \quad \|Q\|_{1} = h \sum_{i\,=\,0}^{N_{x}} |Q_{i}| \, .
\end{equation}
Introducing the sets $\mathcal{L}$ and $\mathcal{K}$ defined by
\begin{equation}
\mathcal{L} = \Big\{ q \in L^{1}\big([0,T) \times \T^{1} \big) \, : \, \int_{\T^{1}} q \, dx = 0 \, \forall \, t \Big\} \, , \quad \mathcal{K} = \big\{ q \in \mathcal{L} \, : \, TV_{T}(q) \leq R \big\} \, ,
\end{equation}
where $R$ is a constant depending on $q_{0}$, it is well known that $\mathcal{K}$ is a compact subset of $\mathcal{L}$. \\

\indent Since the numerical flux involved in (\ref{Roe_Q}) is continuous and satisfies
\begin{equation} \label{consistence}
\mathcal{F}(q,q) = f(q) \quad \forall \, q \, ,
\end{equation}
the considered scheme is consistent with the conservation law (\ref{eq_q}). Hence, applying LeVeque\cite{Finite_volume}, proving that $Q_{h}$ converges to the solution $q$ of (\ref{eq_q}) in $L^{1}\big([0,T) \times \T^{1}\big)$ norm for $h \to 0$ reduces to prove that the scheme is TV-stable. In other words, we need to prove that for $h \in [0,h_{0}]$, the approximation $Q_{h}$ lies in some fixed set $\mathcal{K}$ where $R$ only depends on the initial data $q_{0}$, the final time $T$ and the function $f$. \\

\indent The TV-stability is the consequence of two lemmas.

\begin{lemma} \label{TVD}
The numerical method (\ref{Roe_Q}) is TVD (Total Variation Diminishing), i.e.
\begin{equation}
TV(Q^{n+1}) \leq TV(Q^{n}) \quad \forall \, n \, .
\end{equation}
\end{lemma}

\indent \textit{Proof of the lemma \ref{TVD}}: we develop $Q_{i+1}^{n+1} - Q_{i}^{n+1}$:
\begin{displaymath}
\begin{split}
Q_{i+1}^{n+1} - Q_{i}^{n+1} &= Q_{i+1}^{n} - Q_{i}^{n} - \frac{k}{2h} \Big( f(Q_{i+2}^{n}) - f(Q_{i}^{n}) - f(Q_{i+1}^{n}) + f(Q_{i-1}^{n}) \\
&\qquad \qquad \qquad \qquad \qquad - \big| \alpha(Q_{i+2}^{n}+Q_{i+1}^{n})+\beta \big| (Q_{i+2}^{n}-Q_{i+1}^{n}) \\
&\qquad \qquad \qquad \qquad \qquad + 2 \big| \alpha(Q_{i+1}^{n}+Q_{i}^{n})+\beta \big| (Q_{i+1}^{n}-Q_{i}^{n}) \\
&\qquad \qquad \qquad \qquad \qquad - \big| \alpha(Q_{i}^{n}+Q_{i-1}^{n})+\beta \big| (Q_{i}^{n}-Q_{i-1}^{n}) \Big) \\
&= - \frac{k}{2h} \Big( 1 - \textnormal{sg}\big( \alpha(Q_{i+2}^{n}+Q_{i+1}^{n})+\beta \big) \Big) \big( \alpha(Q_{i+2}^{n}+Q_{i+1}^{n})+\beta \big) (Q_{i+2}^{n} - Q_{i+1}^{n}) \\
& \qquad + \Big( 1 - \frac{k}{h} \, \big| \alpha(Q_{i+1}^{n}+Q_{i}^{n})+\beta \big| \Big) (Q_{i+1}^{n}-Q_{i}^{n}) \\
& \qquad + \frac{k}{2h} \Big( 1 + \textnormal{sg}\big( \alpha(Q_{i}^{n}+Q_{i-1}^{n})+\beta \big) \Big) \big( \alpha(Q_{i}^{n}+Q_{i-1}^{n})+\beta \big) (Q_{i}^{n} - Q_{i-1}^{n}) \, ,
\end{split}
\end{displaymath}
where sg stands for the usual sign function. Using the CFL condition (\ref{CFL}), we write:
\begin{displaymath}
\begin{split}
\big| Q_{i+1}^{n+1} - Q_{i}^{n+1} \big| &\leq \frac{k}{2h} \Big( 1 - \textnormal{sg}\big( \alpha(Q_{i+2}^{n}+Q_{i+1}^{n})+\beta \big) \Big) \big| \alpha(Q_{i+2}^{n}+Q_{i+1}^{n})+\beta \big| \, |Q_{i+2}^{n} - Q_{i+1}^{n}| \\
& \qquad + \Big( 1 - \frac{k}{h} \, \big| \alpha(Q_{i+1}^{n}+Q_{i}^{n})+\beta \big| \Big) \, |Q_{i+1}^{n}-Q_{i}^{n}| \\
& \qquad + \frac{k}{2h} \Big( 1 + \textnormal{sg}\big( \alpha(Q_{i}^{n}+Q_{i-1}^{n})+\beta \big) \Big) \, \big| \alpha(Q_{i}^{n}+Q_{i-1}^{n})+\beta \big| \, |Q_{i}^{n} - Q_{i-1}^{n}| \, .
\end{split}
\end{displaymath}
We deduce then 
\begin{equation}
TV(Q^{n+1}) = \sum_{i\,=\,0}^{N_{x}} \big| Q_{i+1}^{n+1} - Q_{i}^{n+1} \big| \leq \sum_{i\,=\,0}^{N_{x}} \big| Q_{i+1}^{n} - Q_{i}^{n} \big| = TV(Q^{n})
\end{equation}
giving the lemma. $\square$ \\

As a consequence of this lemma, we deduce that there exists a constant $M_{1} > 0$ which only depends on the initial data $q_{0}$ and such that
\begin{equation}
TV(Q^{n}) \leq M_{1} \quad \forall \, n \, .
\end{equation}

\begin{lemma} \label{TV}
For any $n$, we have $\big\|Q^{n+1} - Q^{n} \big\|_{1} \leq 2M_{1}h$. \\
\end{lemma}

\indent \textit{Proof of lemma \ref{TV}:} we develop $Q_{i}^{n+1} - Q_{i}^{n}$:
\begin{displaymath}
\begin{split}
Q_{i}^{n+1} - Q_{i}^{n} &= - \frac{k}{2h} \Big( f(Q_{i+1}^{n}) - f(Q_{i}^{n}) - \big| \alpha(Q_{i+1}^{n}+Q_{i}^{n}) + \beta \big| (Q_{i+1}^{n}-Q_{i}^{n}) \\
&\qquad \qquad \qquad - f(Q_{i}^{n}) + f(Q_{i-1}^{n}) + \big| \alpha(Q_{i}^{n}+Q_{i-1}^{n}) + \beta \big| (Q_{i}^{n}-Q_{i-1}^{n}) \Big) \\
&= -\frac{k}{2h} \Big( 1 - \textnormal{sg}\big( \alpha(Q_{i+1}^{n}+Q_{i}^{n}) + \beta \big) \Big) \big( \alpha(Q_{i+1}^{n}+Q_{i}^{n}) + \beta \big) (Q_{i+1}^{n}-Q_{i}^{n}) \\
&\qquad + \frac{k}{2h} \Big( 1 - \textnormal{sg}\big( \alpha(Q_{i}^{n}+Q_{i-1}^{n}) + \beta \big) \Big) \big( \alpha(Q_{i}^{n}+Q_{i-1}^{n}) + \beta \big) (Q_{i}^{n}-Q_{i-1}^{n}) \, .
\end{split}
\end{displaymath}
Then, using the CFL condition, we obtain
\begin{displaymath}
\begin{split}
\big\|Q^{n+1} - Q^{n} \big\|_{1} &= h \sum_{i\,=\,0}^{N_{x}} \big|Q_{i}^{n+1} - Q_{i}^{n}\big| \\
&\leq h \sum_{i\,=\,0}^{N_{x}} \Bigg[ \frac{k}{2h} \big| \alpha(Q_{i+1}^{n}+Q_{i}^{n}) + \beta \big| \, |Q_{i+1}^{n}-Q_{i}^{n}| \Bigg] \\
&\qquad + h \sum_{i\,=\,0}^{N_{x}} \Bigg[ \frac{k}{2h} \big| \alpha(Q_{i}^{n}+Q_{i-1}^{n}) + \beta \big| \, |Q_{i}^{n}-Q_{i-1}^{n}| \Bigg] \\
&\leq 2 h \sum_{i\,=\,0}^{N_{x}} |Q_{i}^{n}-Q_{i-1}^{n}|
\end{split}
\end{displaymath}
and we complete the proof by applying lemma \ref{TVD}. $\square$ \\ \\

\indent \textit{End of the proof of theorem 5:} combining all these results, we finally have the inequality
\begin{equation}
TV_{T}(Q_{h}) \leq \Big( 1 + \frac{2h}{k} \Big) TM_{1} \, .
\end{equation}
Remarking that $\frac{h}{k}$ is bounded by a constant $M_{2}$ which only depends on the initial data $q_{0}$, we finally obtain that there exists a constant $C > 0$ which only depends on $q_{0}$ and $T$, and such that
\begin{equation}
TV_{T}(Q_{h}) \leq C
\end{equation}
for any space step $h$. We conclude that the numerical method (\ref{Roe_Q}) is TV-stable, yielding the convergence of $Q_{h}$ to $q$. \\

\indent Now, we have to prove that, assuming that $s \geq 3$, the local truncation error is first order accurate. With the notation above, the expression of the local truncation error reads
\begin{equation} \label{LTE}
\begin{split}
e(x_{i},t_{n}) &= \frac{q_{i}^{n} - q_{i}^{n+1}}{k} - \frac{1}{2h} \Big( f(q_{i+1}^{n}) - f(q_{i-1}^{n}) - \big|\alpha(q_{i+1}^{n}+q_{i}^{n})+\beta \big| (q_{i+1}^{n}-q_{i}^{n}) \\ 
&\qquad \qquad \qquad \qquad \qquad \qquad \qquad \qquad + \big| \alpha(q_{i}^{n}+q_{i-1}^{n}) + \beta \big| (q_{i}^{n}-q_{i-1}^{n}) \Big) \, ,
\end{split}
\end{equation}
where $q_{i}^{n} = q(x_{i},t_{n})$. According to the signs of $\alpha(q_{i}^{n}+q_{i-1}^{n}) + \beta$ and $\alpha(q_{i+1}^{n}+q_{i}^{n}) + \beta$, we distinguish four different cases.

\begin{itemize}
\item First case: we assume that $\alpha(q_{i+1}^{n}+q_{i}^{n}) + \beta \geq 0$ and $\alpha(q_{i}^{n}+q_{i-1}^{n}) + \beta \geq 0$. Then, the local truncation error becomes
\begin{equation} \label{case1}
e(x_{i},t_{n}) = \frac{q_{i}^{n} - q_{i}^{n+1}}{k} - \frac{f(q_{i}^{n}) - f(q_{i-1}^{n})}{h} \, .
\end{equation}
Since $s \geq 3$, $q \in \Cl^{2}\big([0,T) \times \T^{1}\big)$, so we can write
\begin{equation}
e(x_{i},t_{n}) = -\D_{t} q(x_{i},t_{n}) - \D_{x} \Big( f\big(q(x_{i},t_{n}) \big) \Big) + \mathcal{O}(k) + \mathcal{O}(h) = \mathcal{O}(h+k) \, .
\end{equation}

\item Second case: we assume that $\alpha(q_{i+1}^{n}+q_{i}^{n}) + \beta \leq 0$ and $\alpha(q_{i}^{n}+q_{i-1}^{n}) + \beta \leq 0$. Then, the local truncation error becomes
\begin{equation} \label{case2}
e(x_{i},t_{n}) = \frac{q_{i}^{n} - q_{i}^{n+1}}{k} - \frac{f(q_{i+1}^{n})-f(q_{i}^{n})}{h} \, .
\end{equation}
Since $s \geq 3$, $q \in \Cl^{2}\big([0,T) \times \T^{1}\big)$, so we can write
\begin{equation}
e(x_{i},t_{n}) = -\D_{t} q(x_{i},t_{n}) - \D_{x} \Big( f\big(q(x_{i},t_{n}) \big) \Big) + \mathcal{O}(k) + \mathcal{O}(h) = \mathcal{O}(h+k) \, .
\end{equation}

\item Third case: we assume that $\alpha(q_{i+1}^{n}+q_{i}^{n}) + \beta \geq 0$ and $\alpha(q_{i}^{n}+q_{i-1}^{n}) + \beta \leq 0$. Then, the function $\lambda$ defined by
\begin{equation} \label{lambda}
\lambda(x) = \alpha \big( q(x,t_{n}) + q(x-h,t_{n}) \big) + \beta
\end{equation}
admits a zero denoted $x_{*}$ in $[x_{i},x_{i+1}]$. Since $s \geq 3$, $\lambda$ is of class $\Cl^{2}$ and verifies
\begin{equation}
\begin{split}
\lambda(x_{i}) &= \lambda(x_{*}) + \mathcal{O}(x_{i}-x_{*}) = \mathcal{O}(h) \, , \\ \lambda(x_{i+1}) &= \lambda(x_{*}) + \mathcal{O}(x_{i+1}-x_{*}) = \mathcal{O}(h) \, .
\end{split}
\end{equation}
In the same way, we have
\begin{equation} \label{lambda_0}
f'(q_{i}^{n}) = 2\alpha q_{i}^{n} + \beta = \lambda(x_{i+1}) + \mathcal{O}(h) = \mathcal{O}(h) \, .
\end{equation}
The local truncation error is of the form
\begin{equation}
e(x_{i},t_{n}) = \frac{q_{i}^{n}-q_{i}^{n+1}}{k} = -\D_{t}q(x_{i},t_{n}) + \mathcal{O}(k) = f'\big(q(x_{i},t_{n})\big)\, \D_{x}q(x_{i},t_{n}) + \mathcal{O}(k) \, .
\end{equation}
Then, using (\ref{lambda_0}), we obtain
\begin{equation} \label{case3}
e(x_{i},t_{n}) = \mathcal{O}(h+k) \, .
\end{equation}
\item Fourth case: we assume that $\alpha(q_{i+1}^{n}+q_{i}^{n}) + \beta \leq 0$ and $\alpha(q_{i}^{n}+q_{i-1}^{n}) + \beta \geq 0$. We proceed as we did in the precedent case, i.e. we deduce that the function $\lambda$ defined by (\ref{lambda}) admits a zero $x_{*} \in [x_{i},x_{i+1}]$, and then the result (\ref{lambda_0}) is true. Hence,
\begin{equation} \label{case4}
\begin{split}
e(x_{i},t_{n}) &= \frac{q_{i}^{n}-q_{i}^{n+1}}{k} - \frac{f(q_{i+1}^{n}) - f(q_{i-1}^{n})}{h} \\
&= -\D_{t}q(x_{i},t_{n}) - 2 \D_{x}\Big(f\big(q(x_{i},t_{n})\big)\Big) + \mathcal{O}(h+k) \\
&= -f'(q_{i}^{n})\D_{x}q(x_{i},t_{n}) + \mathcal{O}(h+k) \\
&= \mathcal{O}(h+k) \, .
\end{split}
\end{equation}
\end{itemize}

We conclude that, in any case, the local truncation error is $\mathcal{O}(h+k)$, i.e. first order accurate. $\square$

\section{Numerical results}

\setcounter{equation}{0}

\indent The first goal of this section is to numerically show that $u^{\epsilon}(x,t)-U_{h}(x,\frac{t}{\epsilon},t) = \mathcal{O}(\epsilon)$ and $\rho^{\epsilon}(x,t) - R_{h}(x,\frac{t}{\epsilon},t) = \mathcal{O}(\epsilon)$. Secondly, we will briefly analyze the gain in term of CPU time our two-scale numerical method brings when compared with a classical method consisting in solving directly (\ref{Grenier}). Lastly, by using our method, we will explore the simulation of experiments which were not accessible before because of too small Mach number.

\subsection{Convergence in $\epsilon$}

\indent We consider here the initial data $u_{0}$ and $\rho_{0}$ defined by
\begin{equation}\label{init}
u_{0}(x) = 1 + \frac{\cos(x)}{2} \quad \textnormal{and} \quad \rho_{0}(x) = 1 + \frac{\sin(x)}{2} \, ,
\end{equation}
and we take a uniform mesh with $N_{x}+1$ points denoted $x_{0},\dots,x_{N_{x}}$. \\
\indent As illustrated in figures 1 and 2, the functions $U_{h}(x,\frac{t}{\epsilon},t)$ and $R_{h}(x,\frac{t}{\epsilon},t)$ are very close to $u^{\epsilon}$ and $\rho^{\epsilon}$ respectively. The numerical experiment showed in these figures is made with $\epsilon = 0.05$, $\gamma = 1$, and $N_{x} = 1023$, $T = 2.5$.

\begin{center}
\includegraphics[scale=0.28,angle=270]{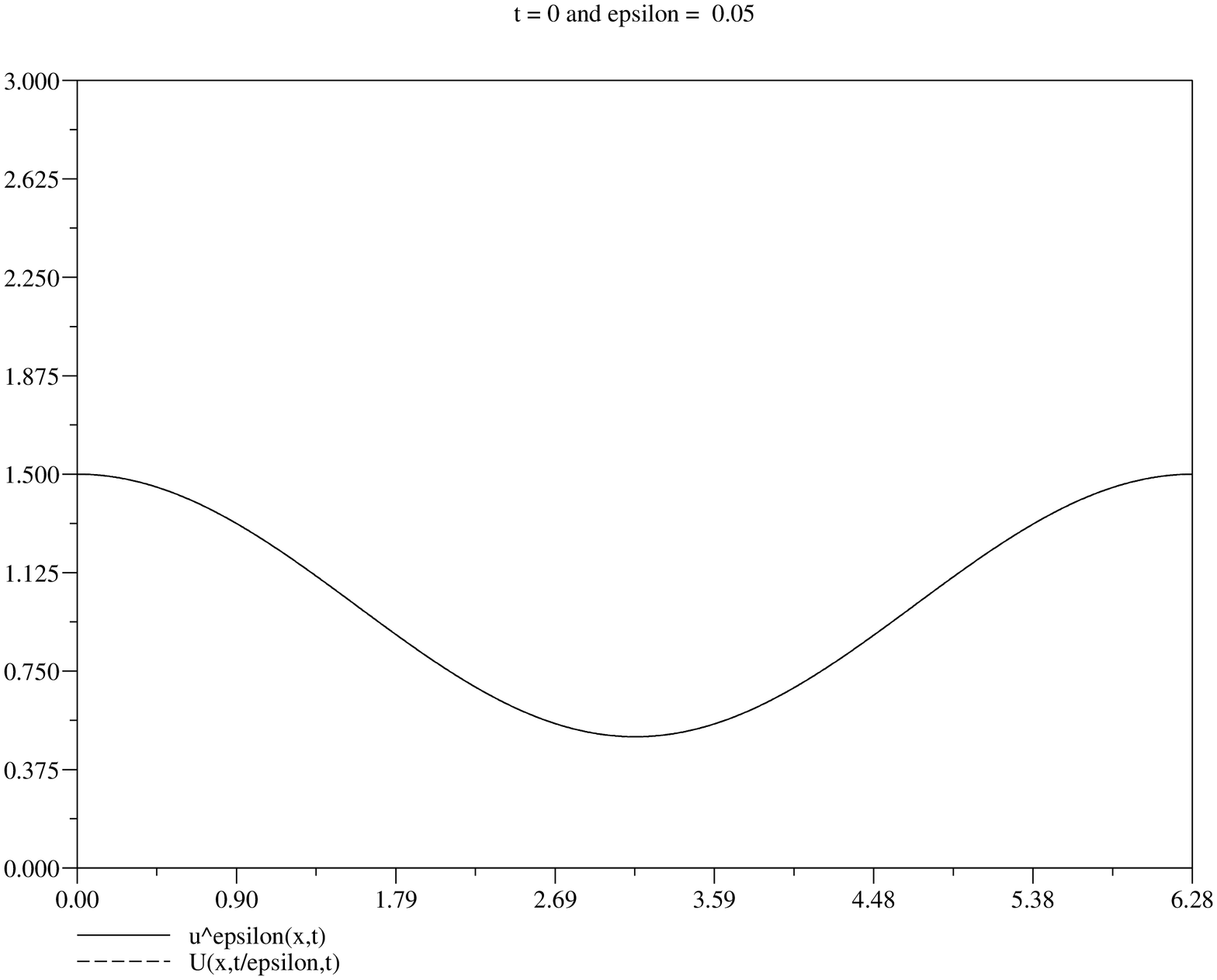} \hspace{-1.15cm}
\includegraphics[scale=0.28,angle=270]{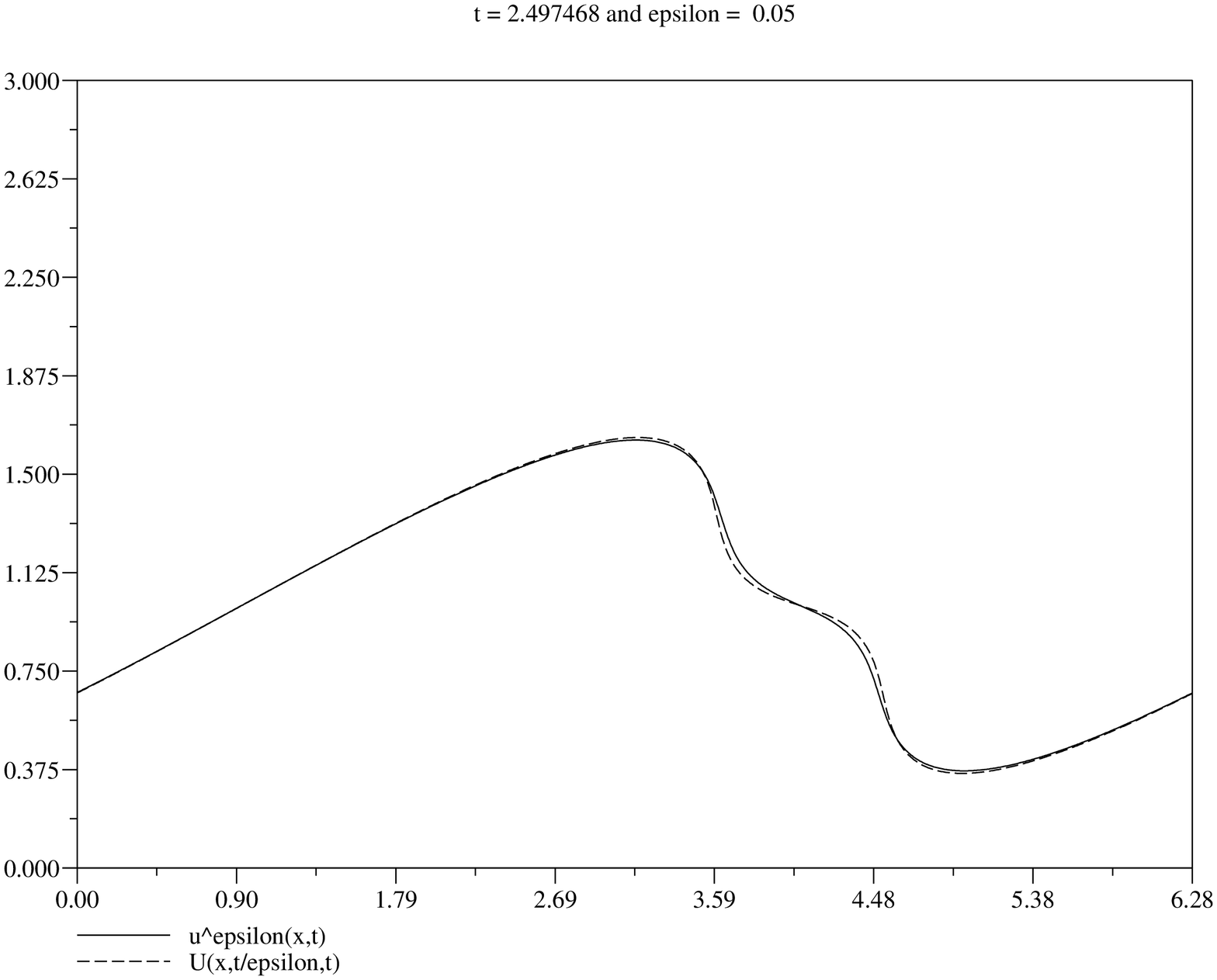} \\
\caption{Comparison of $u^{\epsilon}(\cdot,t)$ and $U_{h}(\cdot,\frac{t}{\epsilon},t)$ at times $t = 0$ and $t = 2.5$.}
\includegraphics[scale=0.28,angle=270]{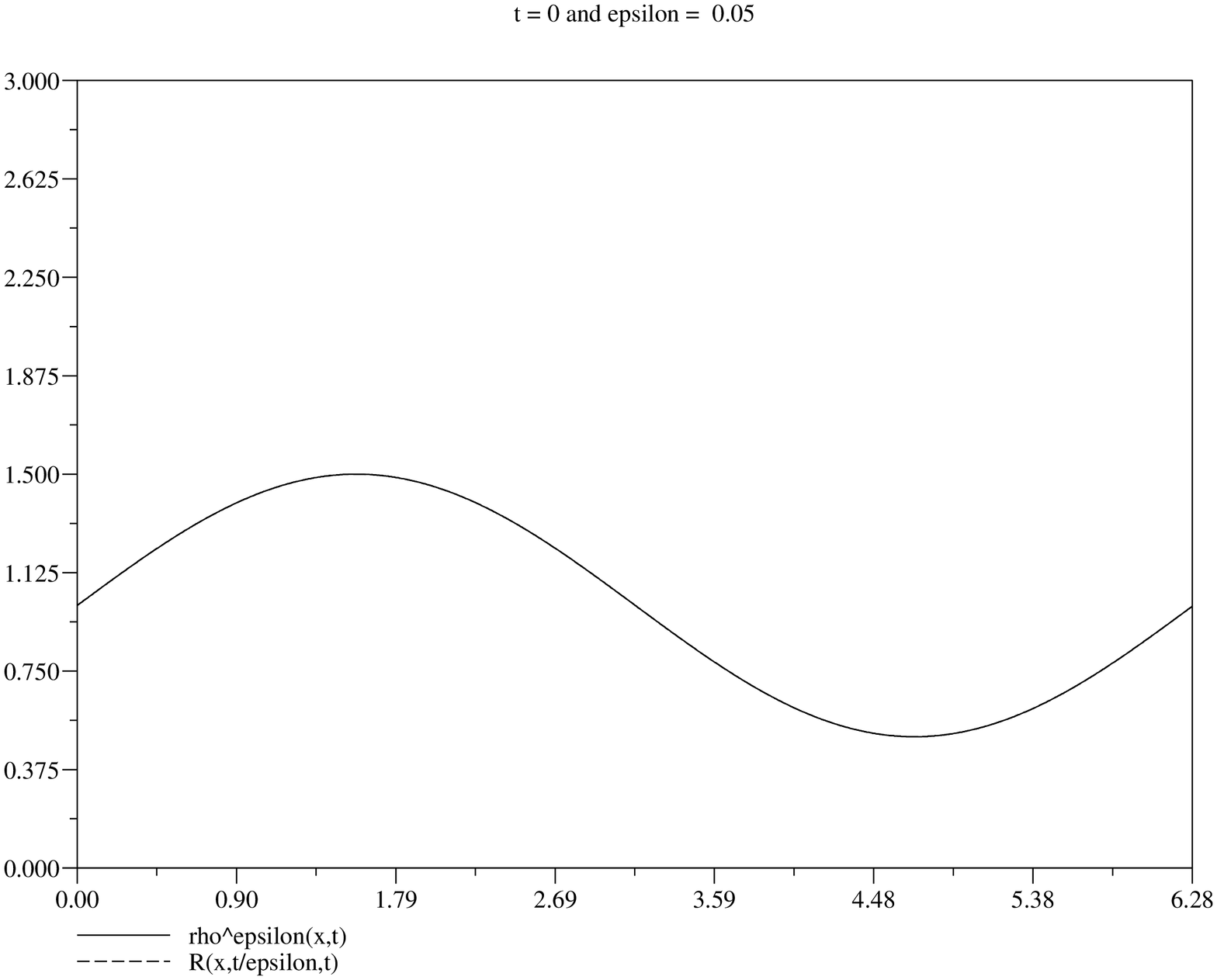} \hspace{-1.15cm}
\includegraphics[scale=0.28,angle=270]{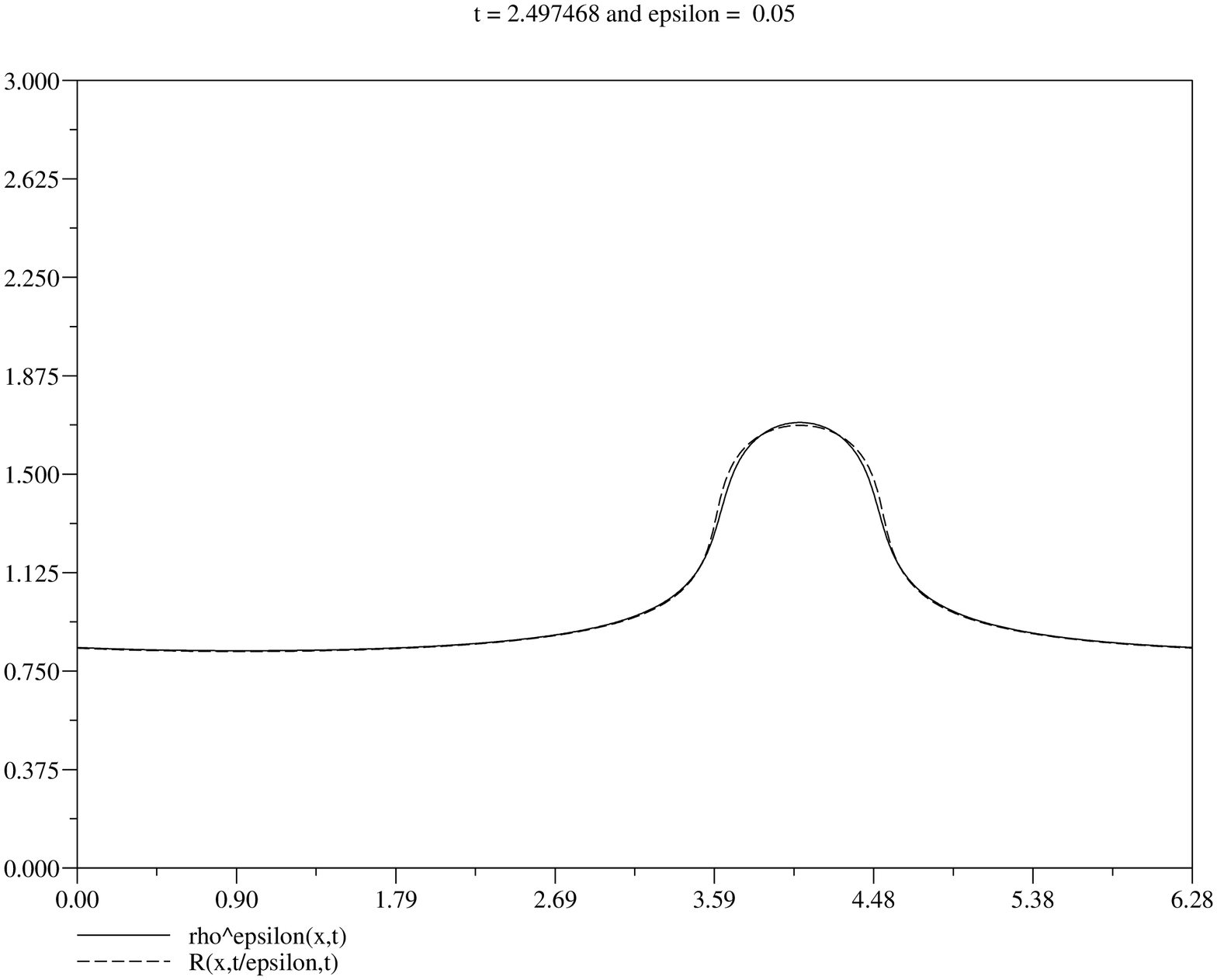} \\
\caption{Comparison of $\rho^{\epsilon}(\cdot,t)$ and $R_{h}(\cdot,\frac{t}{\epsilon},t)$ at times $t = 0$ and $t = 2.5$.}
\end{center}

\indent In order to quantify this good accuracy of $U_{h}(x,\frac{t}{\epsilon},t)$ with $u^{\epsilon}$ and $R_{h}(x,\frac{t}{\epsilon},t)$ with $\rho^{\epsilon}$, we compute the errors $u^{\epsilon}(x,t)-U_{h}(x,\frac{t}{\epsilon},t)$ and $\rho^{\epsilon}(x,t)-R_{h}(x,\frac{t}{\epsilon},t)$ in $L^{p}\big([0,T) \times \T^{1}\big)$ norm ($p = 1, 2, \infty$), for several values of $\epsilon$ ranging in $\{ 0.01, 0.03, 0.05, 0.07, 0.1\}$, and we obtain the figures and the array below:
\vspace{-0.7cm}
\begin{center}
\includegraphics[scale=0.45,angle=270]{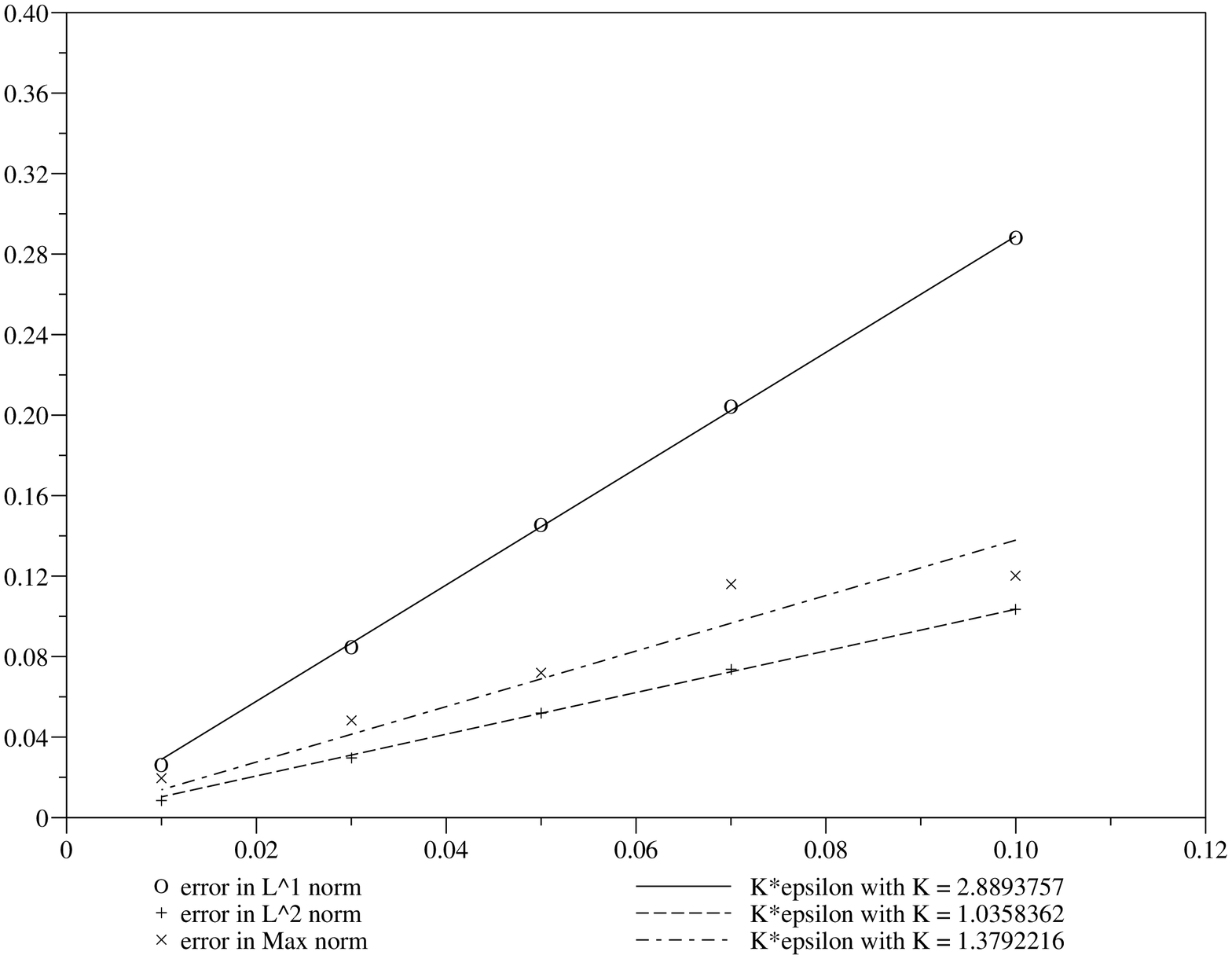}
\caption{Error $u^{\epsilon}(x,t) - U_{h}(x,\frac{t}{\epsilon},t)$ in $L^{1}$, $L^{2}$ and $L^{\infty}$ norms.}
\includegraphics[scale=0.45,angle=270]{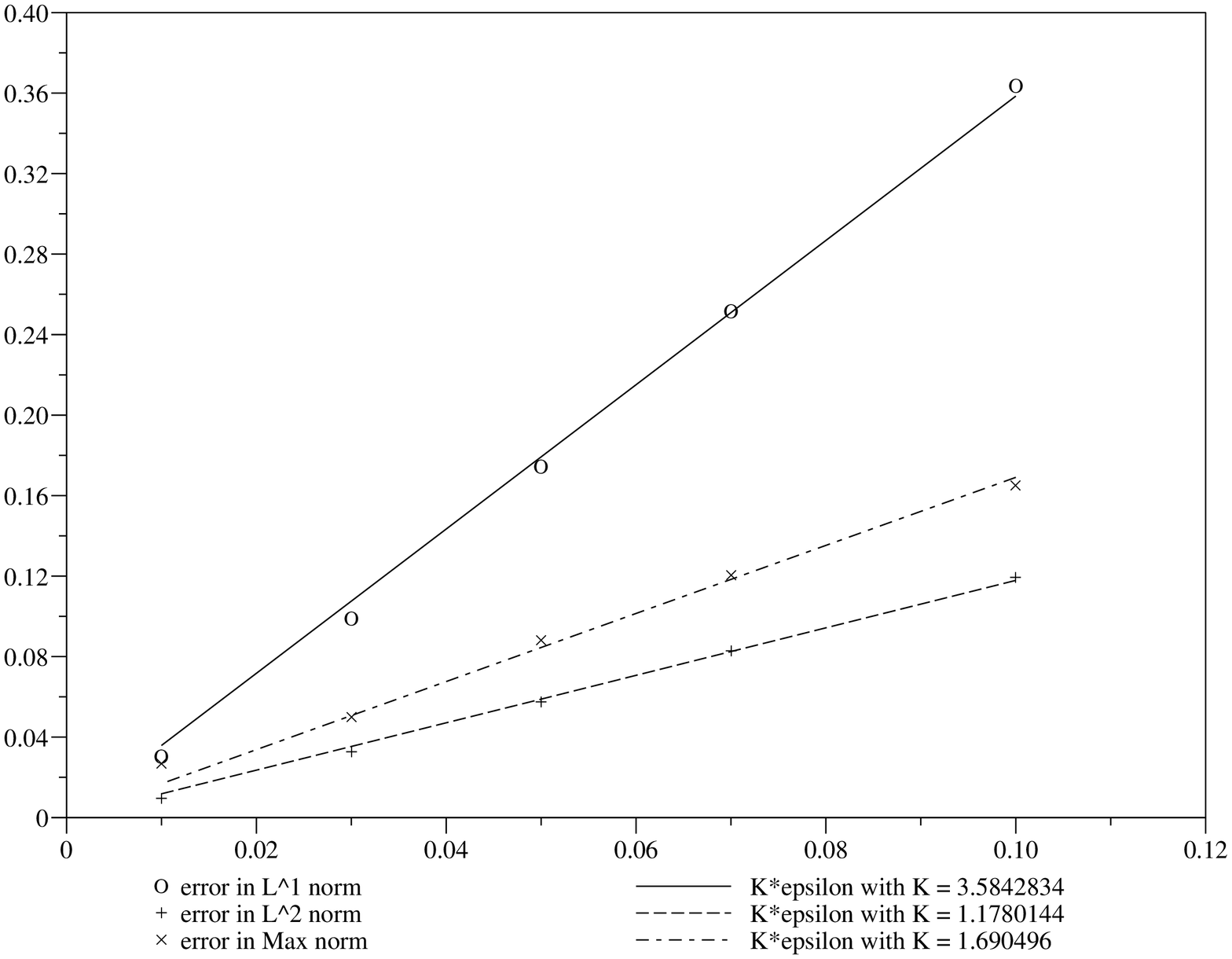}
\caption{Error $\rho^{\epsilon}(x,t) - R_{h}(x,\frac{t}{\epsilon},t)$ in $L^{1}$, $L^{2}$ and $L^{\infty}$ norms.}
\end{center}
\begin{center}
\begin{tabular}{|c|c|c|c|c|c|c|}
\hline
Value & \multicolumn{3}{|c|}{Error $u^{\epsilon}(x,t)-U_{h}(x,\frac{t}{\epsilon},t)$ } & \multicolumn{3}{|c|}{Error $\rho^{\epsilon}(x,t)-R_{h}(x,\frac{t}{\epsilon},t)$ } \\
\cline{2-7}
of $\epsilon$ & $L^{1}$ norm & $L^{2}$ norm & $L^{\infty}$ norm & $L^{1}$ norm & $L^{2}$ norm & $L^{\infty}$ norm \\
\hline
0.1 & 0.2880016 & 0.1034517 & 0.1201132 & 0.3634635 & 0.1192177 & 0.1651374 \\
\hline
0.07 & 0.2043288 & 0.0735498 & 0.1160378 & 0.2516643 & 0.0828920 & 0.1203342 \\
\hline
0.05 & 0.1452756 & 0.0518309 & 0.0719549 & 0.1744198 & 0.0575395 & 0.0879887 \\
\hline
0.03 & 0.0845621 & 0.0296311 & 0.0483544 & 0.0987818 & 0.0326775 & 0.0499608 \\
\hline
0.01 & 0.0260695 & 0.0085251 & 0.0195340 & 0.0303524 & 0.0092954 & 0.0269734 \\
\hline
\end{tabular}
\vspace{-0.5cm}
\begin{table}[ht]
\caption{Errors $u^{\epsilon}(x,t)-U_{h}(x,\frac{t}{\epsilon},t)$ and $\rho^{\epsilon}(x,t)-R_{h}(x,\frac{t}{\epsilon},t)$.}
\end{table}
\end{center}
\vspace{-0.5cm}
\indent The results above show an error $u^{\epsilon}(x,t) - U_{h}(x,\frac{t}{\epsilon},t)$ in $L^{1}$ norm decreasing when $\epsilon  \to 0$ as $K_{1} \epsilon$ with $K_{1} \approx 2.8893757$, as $K_{2}\epsilon$ with $K_{2} \approx 1.0358362$ in $L^{2}\big([0,T) \times \T^{1} \big)$ norm, and and as $K_{\infty}\epsilon$ with $K_{\infty} \approx 1.3792216$ in $L^{\infty}\big([0,T) \times \T^{1} \big)$ norm. Concerning the error $\rho^{\epsilon}(x,t)-R_{h}(x,\frac{t}{\epsilon},t)$, we observe that it decreases when $\epsilon \to 0$ as $K_{1}' \epsilon$ with $K_{1}' \approx 3.5842834$ in $L^{1}\big([0,T) \times \T^{1}\big)$ norm, as $K_{2}'\epsilon$ with $K_{2}' \approx 1.1780144$ in $L^{2}\big([0,T) \times \T^{1} \big)$ norm, and and as $K_{\infty}'\epsilon$ with $K_{\infty}' \approx 1.690496$ in $L^{\infty}\big([0,T) \times \T^{1} \big)$ norm.

\subsection{CPU time cost}

\indent One of the motivations of the present paper is the high CPU time cost of classical methods when used with small Mach number, as we explained it in the introduction. For example, if we apply Roe's finite volume scheme on the model (\ref{Grenier}), we must verify the CFL condition 
\begin{equation}
\frac{k}{h} \max_{i\, , \, n} \textstyle \Big| \hat{u}_{i}^{n} \pm \frac{\sqrt{\overline{P}_{i}^{n}}}{\epsilon} \Big| \leq 1 \, ,
\end{equation}
where $h$ is the space step, $k$ is the time step, $\hat{u}_{i}^{n}$ and $\overline{P}_{i}^{n}$ are defined by
\begin{equation}
\begin{split}
\hat{u}_{i}^{n} &= \frac{u_{i}^{n} \sqrt{1+\epsilon \rho_{i}^{n}} + u_{i-1}^{n} \sqrt{1+\epsilon \rho_{i-1}^{n}}}{\sqrt{1+\epsilon \rho_{i}^{n}} + \sqrt{1+\epsilon \rho_{i-1}^{n}}} \, , \\
\overline{P}_{i}^{n} &= \left\{
\begin{array}{ll}
\frac{(1+\epsilon \rho_{i}^{n})^{\gamma} - (1+\epsilon \rho_{i-1}^{n})^{\gamma}}{\gamma \, \epsilon (\rho_{i}^{n} - \rho_{i-1}^{n})} & \textnormal{if $\rho_{i}^{n} \neq \rho_{i-1}^{n}$} \, , \\
(1+\epsilon \rho_{i}^{n})^{\gamma-1} & \textnormal{else} \, .
\end{array}
\right. 
\end{split}
\end{equation}

Hence, in order to garantee the stability of the scheme, the smaller $\epsilon$ is, the smaller $k$ must be (and the higher the CPU time cost will be). One advantage of our two-scale numerical method is that its CFL condition (\ref{CFL}) does not depend on $\epsilon$, so $k$ does not have to diminish with $\epsilon$. \\
\vspace{-0.3cm}
\begin{center}
\begin{table}[ht]
\begin{tabular}{|c|c|c|c|}
\hline
Value & Roe's method & \multicolumn{2}{|c|}{Two-scale numerical method on (\ref{def_UR})-(\ref{def_FB})} \\
\cline{3-4}
of $\epsilon$ & on (\ref{Grenier}) & Computation of $(F_{h},B_{h})$ & Computation of $(U_{h},R_{h})$ \\
\hline
0.1 & 24 m 34 s 36'' & 3 m 45 s 31'' & 42 s 79'' \\
\hline
0.07 & 32 m 55 s 44'' & 3 m 45 s 22'' & 42 s 93'' \\
\hline
0.05 & 44 m 7 s 21'' & 3 m 45 s 38'' & 42 s 74'' \\
\hline
0.03 & 1 h 10 m 15 s & 3 m 45 s 34'' & 42 s 80'' \\
\hline
0.01 & 3 h 20 m 53 s 78'' & 3 m 45 s 12'' & 42 s 91'' \\
\hline
\end{tabular}
\caption{Comparison in terms of CPU time cost between Roe's method on (\ref{Grenier}) and the two-scale numerical method on (\ref{def_UR})-(\ref{def_FB}) with the initial data (\ref{init}), $N_{x} = 1023$, $\gamma = 1$.}
\end{table}
\end{center}

\newpage

\indent The array above gives the comparison in terms of CPU time cost between Roe's method applied on the non-homogenized model (\ref{Grenier}), and our two-scale numerical method. These results have been obtained with the computer characteristics below:
\vspace{-0.2cm}
\begin{itemize}
\item processor: Pentium$^{\copyright}$ M 715,
\vspace{-0.6cm}
\item memory: 512 Mo DDR-RAM,
\vspace{-0.6cm}
\item operating system: SuSe$^{\copyright}$ 9.1 Pro,
\vspace{-0.6cm}
\item compiler: gcc 3.3.3-41.
\end{itemize}

In particular, we distinguish the time cost for the computation of $F_{h}$ and $B_{h}$ on one hand, and the time cost for the reconstruction of $U_{h}(x,\frac{t}{\epsilon},t)$ and $R_{h}(x,\frac{t}{\epsilon},t)$ on the other hand.

\subsection{Numerical experiments with small Mach number}

In this paragraph, we present some numerical results obtained with our two-scale numerical method on simulations of experiments inducing a very small Mach number. These results are obtained with the initial data (\ref{init}), $N_{x}+1 = 1024$ points in $x$, and $\gamma = 1$.
\vspace{-0.5cm}
\begin{center}
\includegraphics[scale=0.26,angle=270]{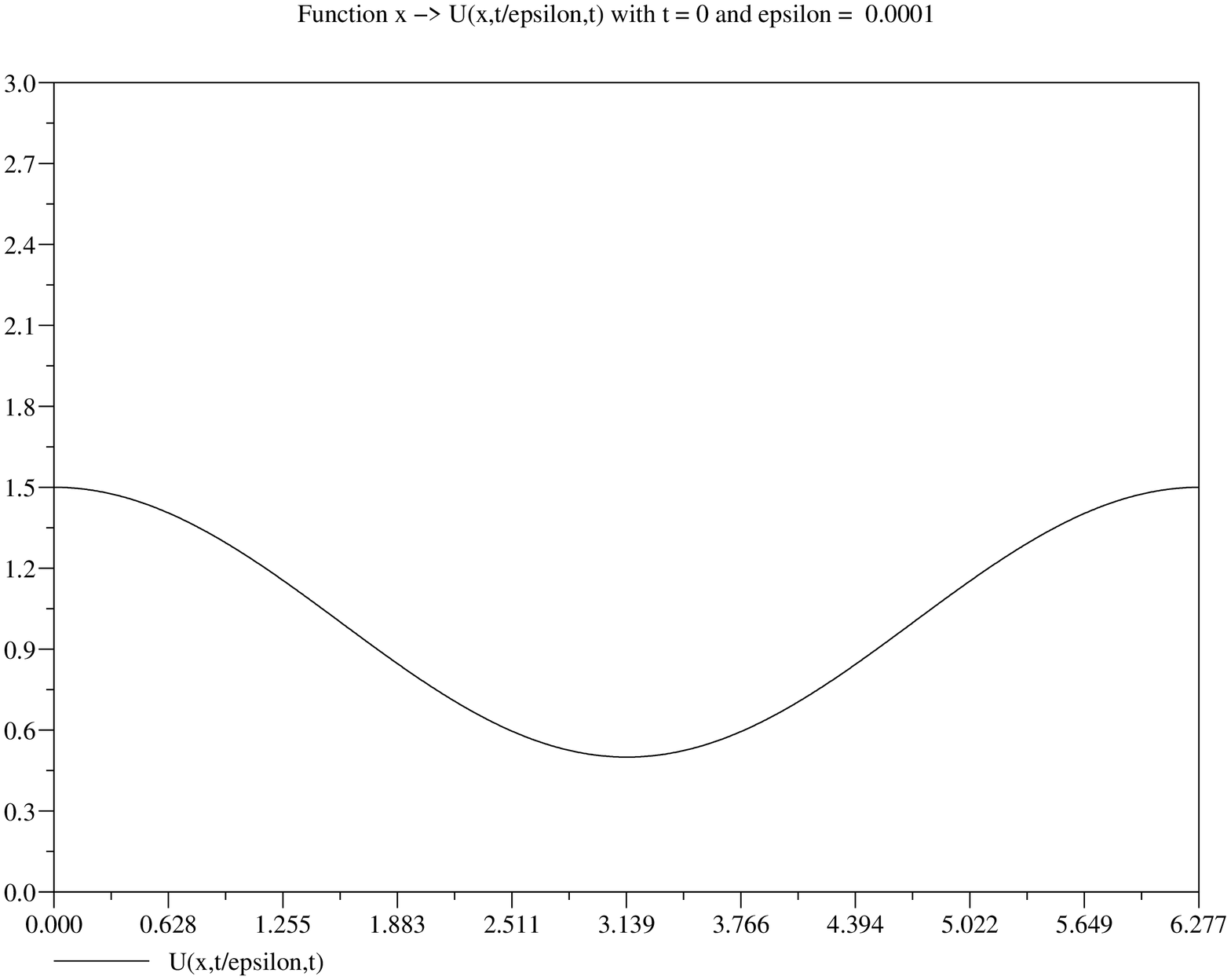} \hspace{-1.15cm}
\includegraphics[scale=0.26,angle=270]{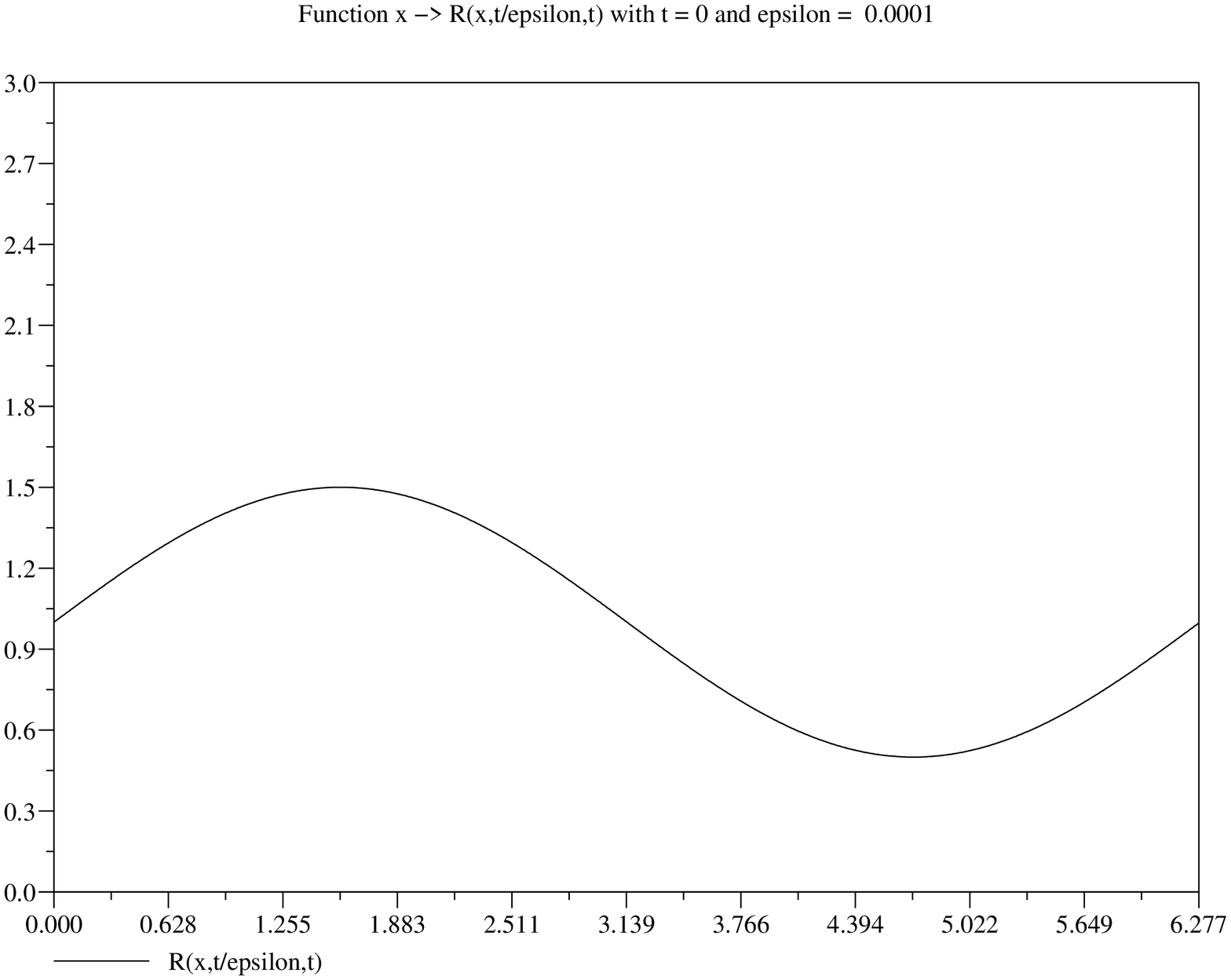} \\
\vspace{-0.5cm}
\includegraphics[scale=0.26,angle=270]{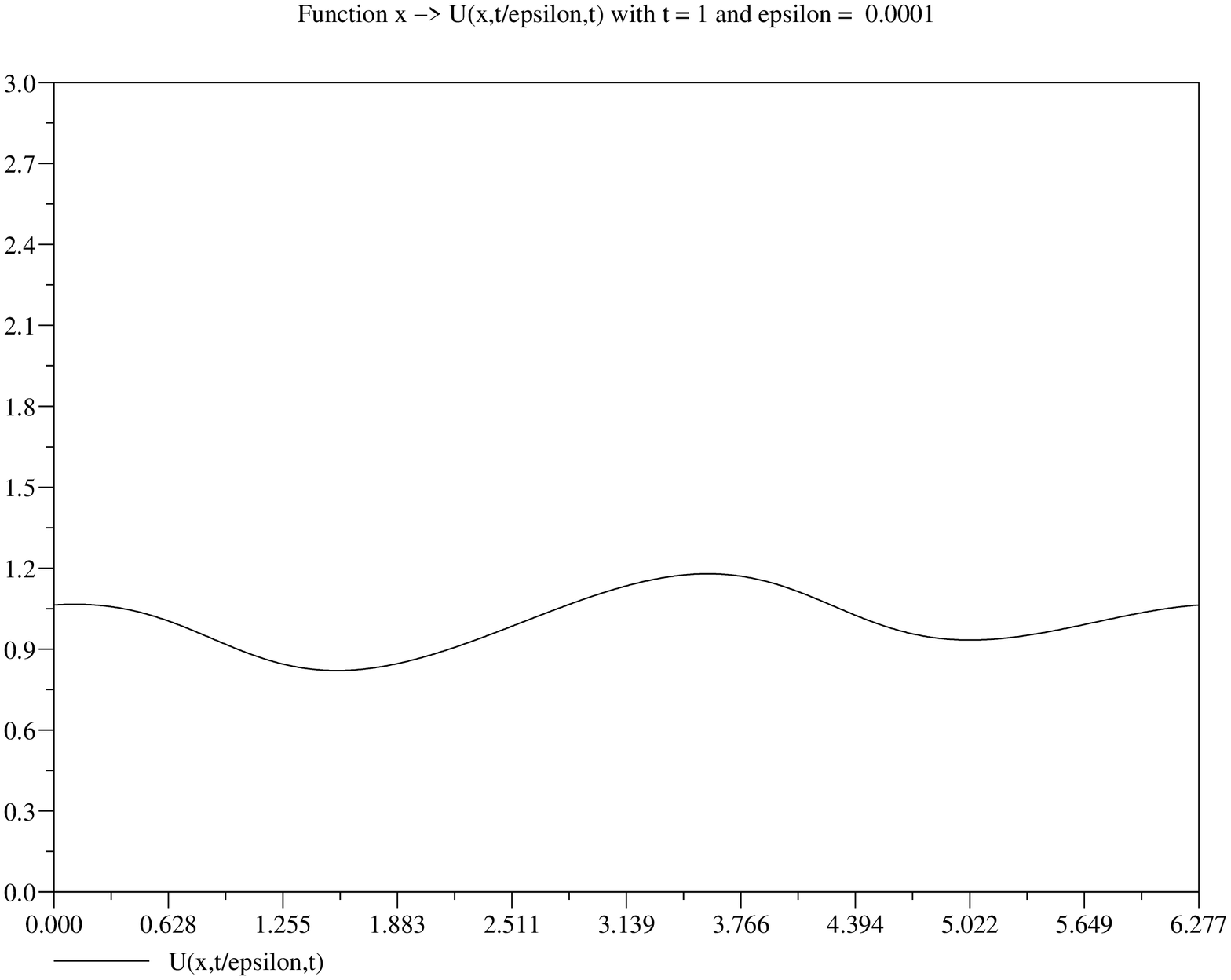} \hspace{-1.15cm}
\includegraphics[scale=0.26,angle=270]{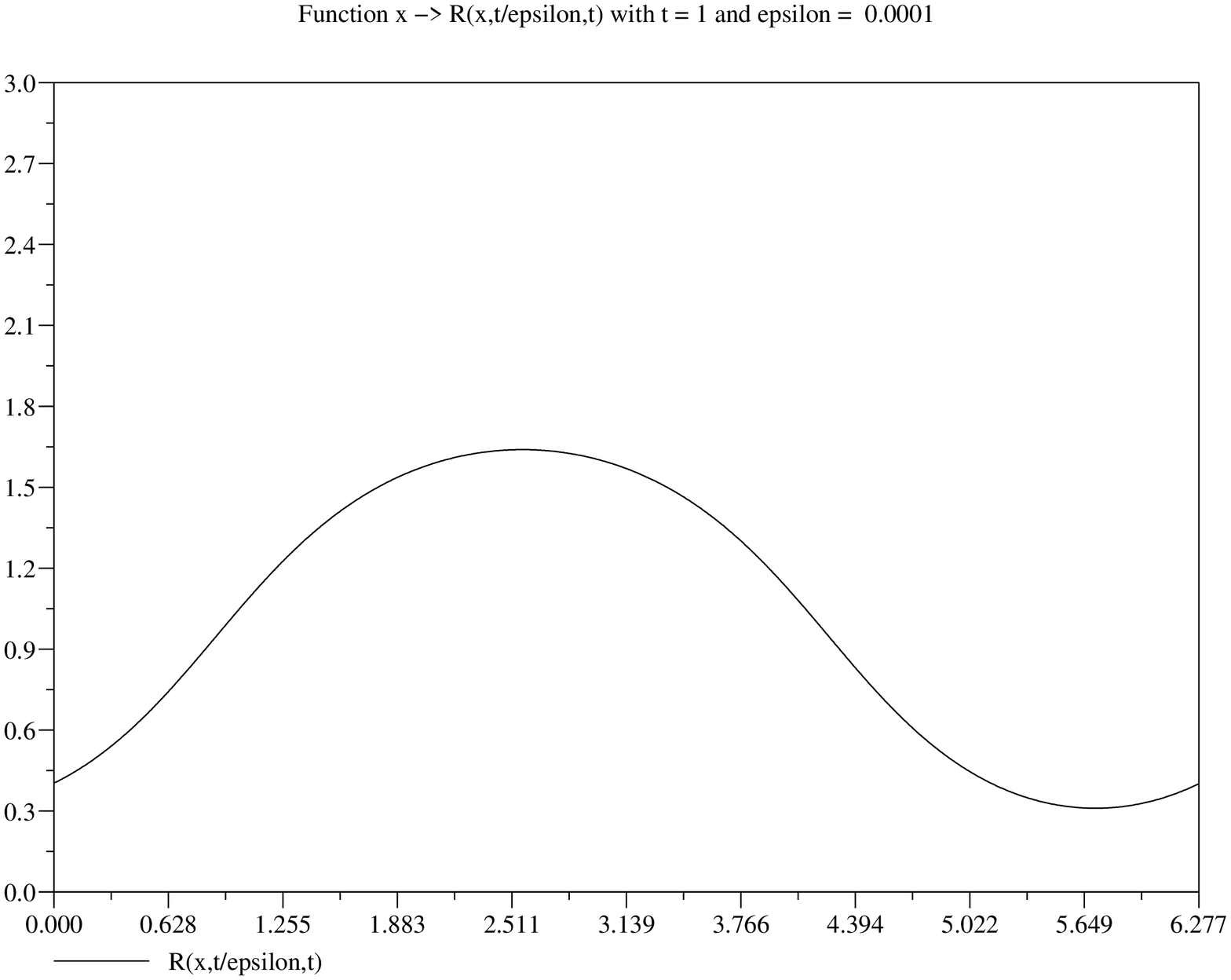}
\vspace{-0.5cm}
\caption{$U_{h}(x,\frac{t}{\epsilon},t)$ (left) and $R_{h}(x,\frac{t}{\epsilon},t)$ (right) with $\epsilon = 10^{-4}$ at times $t = 0$, $t = 1$.}
\includegraphics[scale=0.26,angle=270]{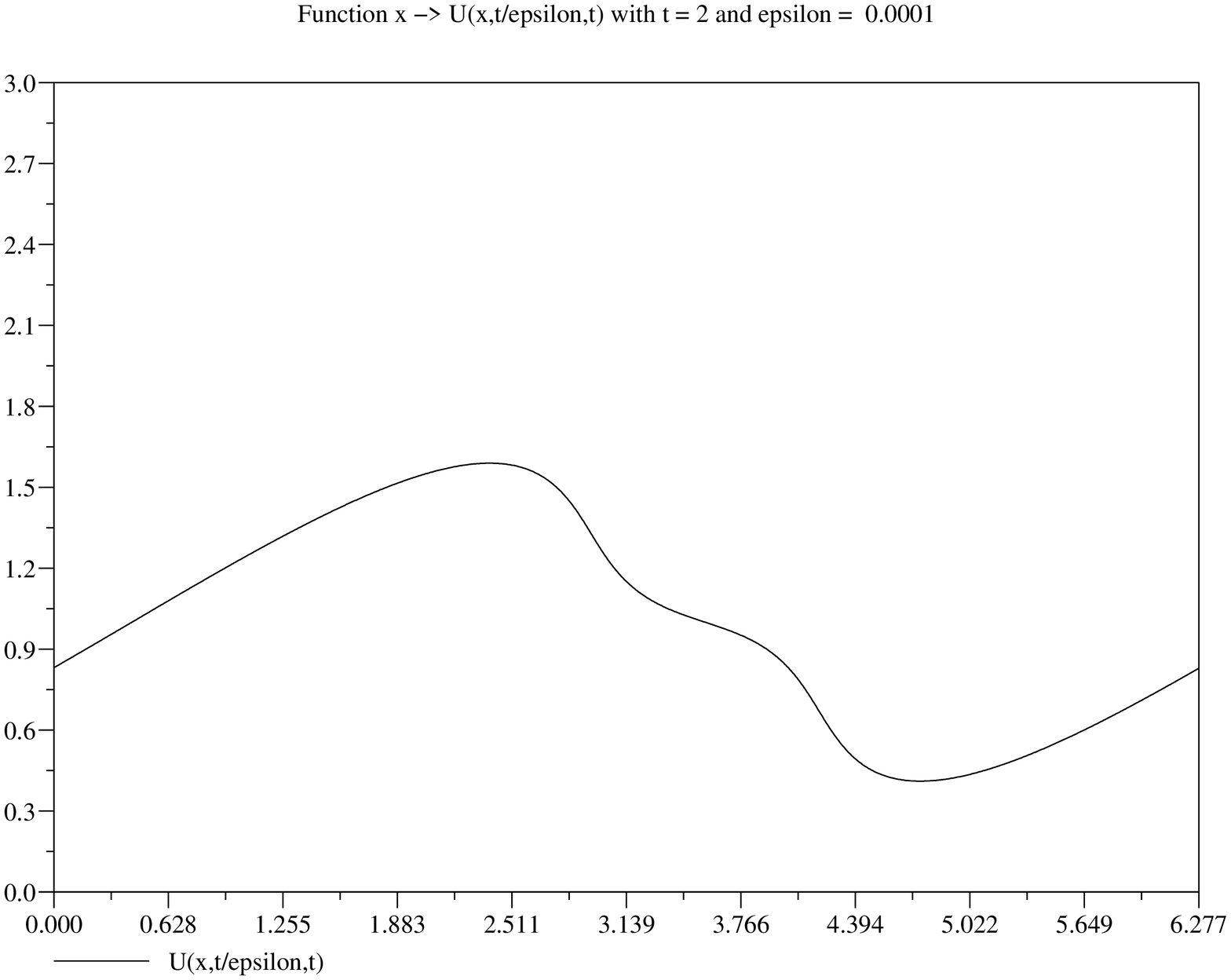} \hspace{-1.15cm}
\includegraphics[scale=0.26,angle=270]{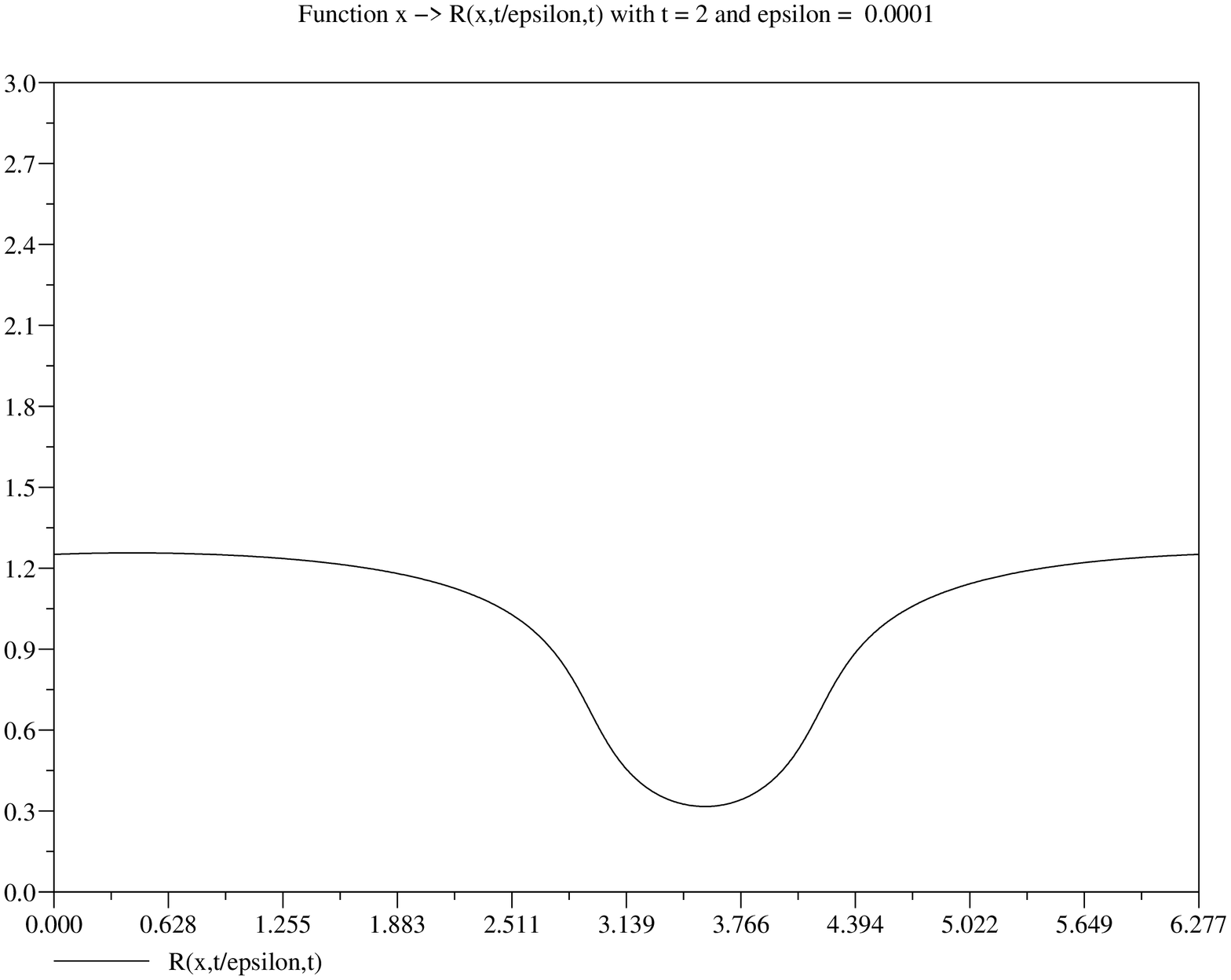}
\vspace{-0.5cm}
\caption{$U_{h}(x,\frac{t}{\epsilon},t)$ (left) and $R_{h}(x,\frac{t}{\epsilon},t)$ (right) with $\epsilon = 10^{-4}$ at time $t = 2$.}
\includegraphics[scale=0.26,angle=270]{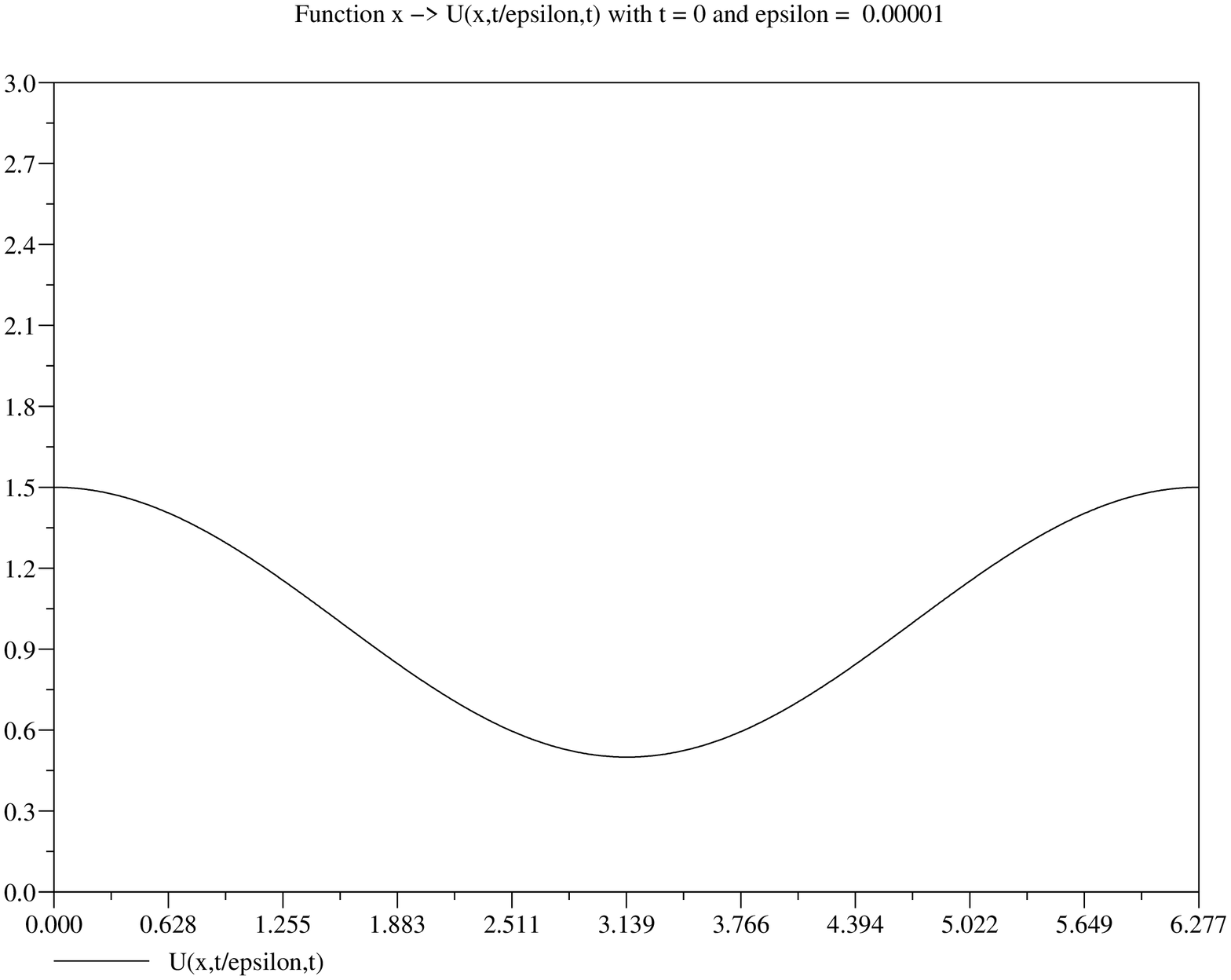} \hspace{-1.15cm} \includegraphics[scale=0.26,angle=270]{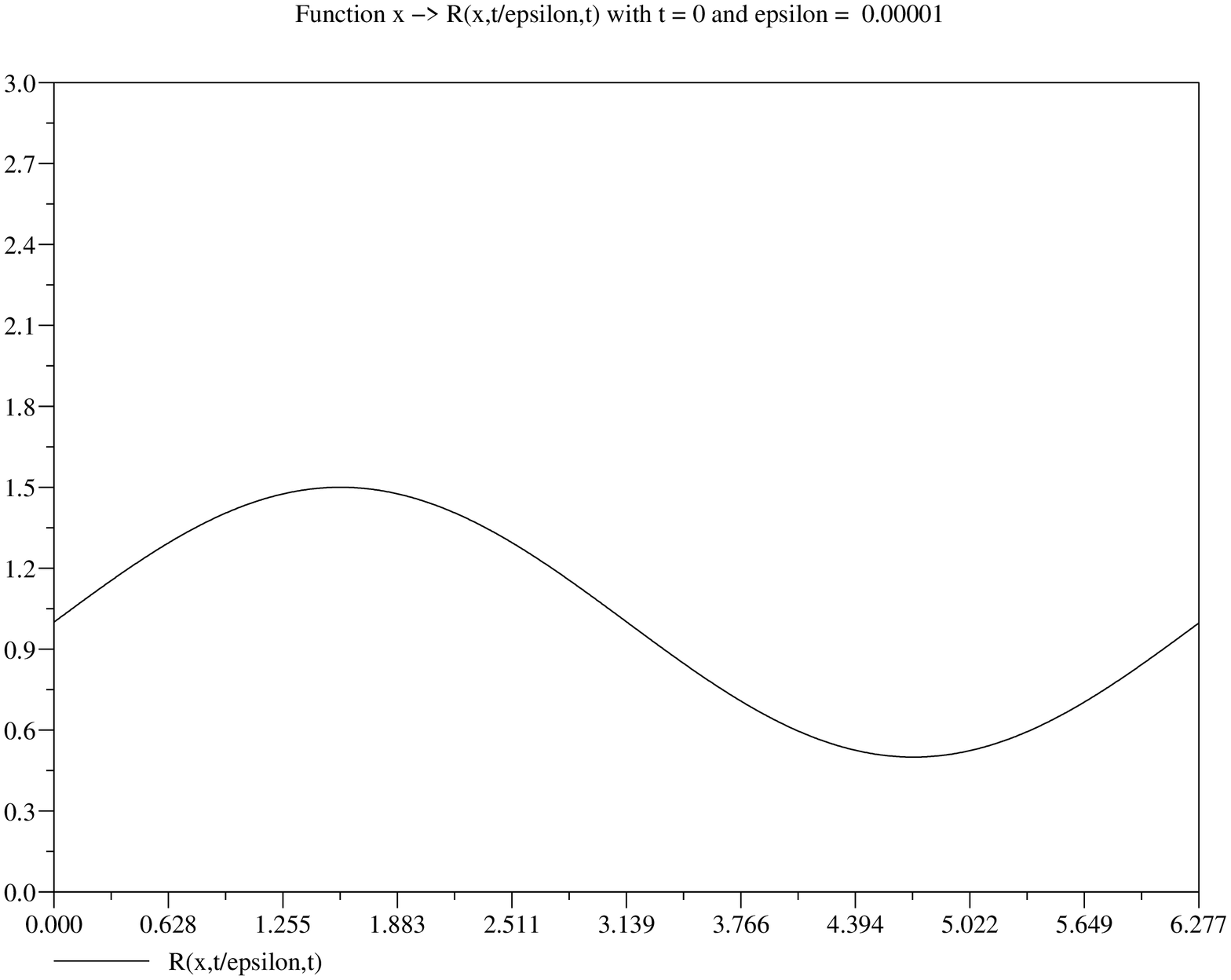} \\
\vspace{-0.5cm}
\includegraphics[scale=0.26,angle=270]{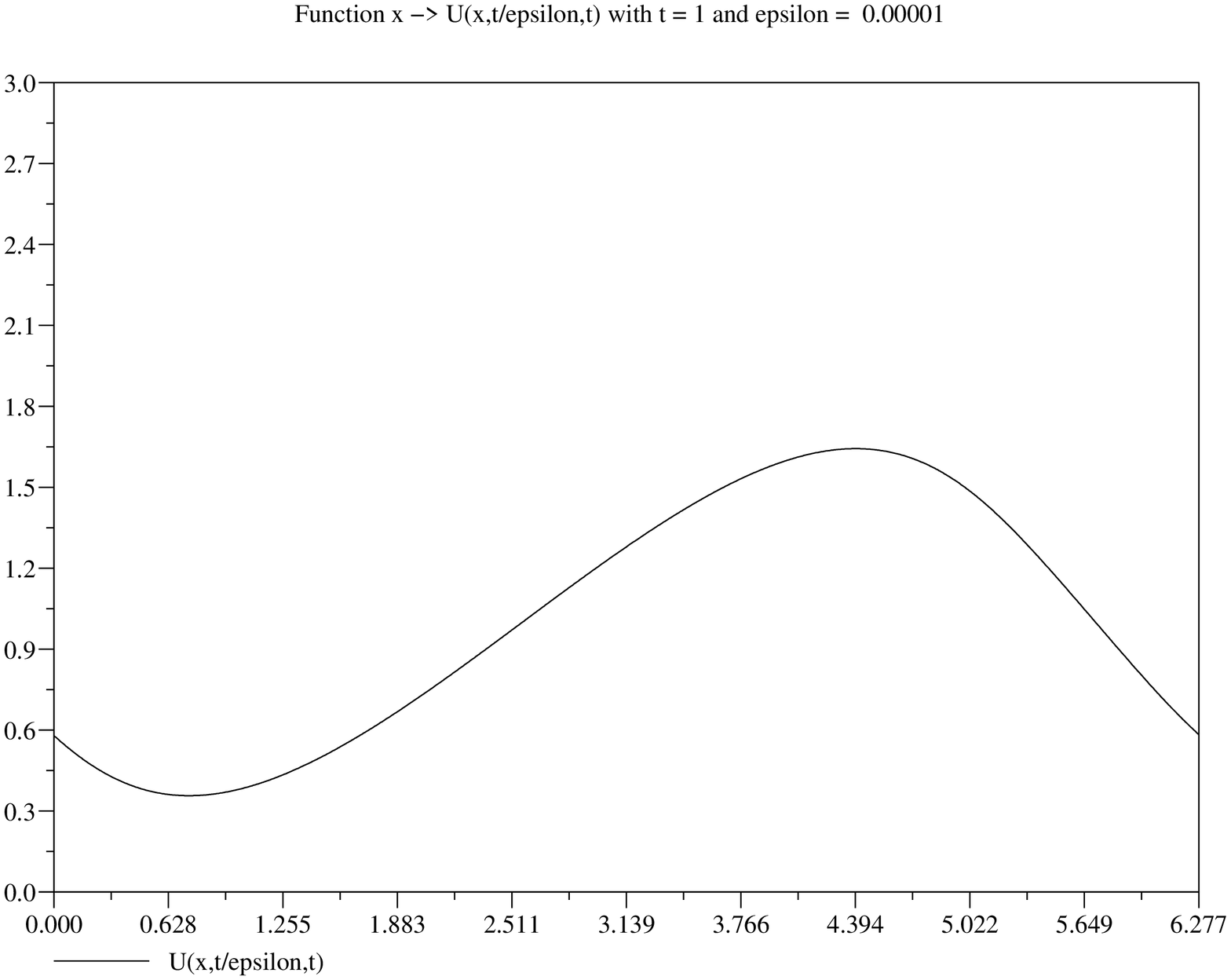} \hspace{-1.15cm}
\includegraphics[scale=0.26,angle=270]{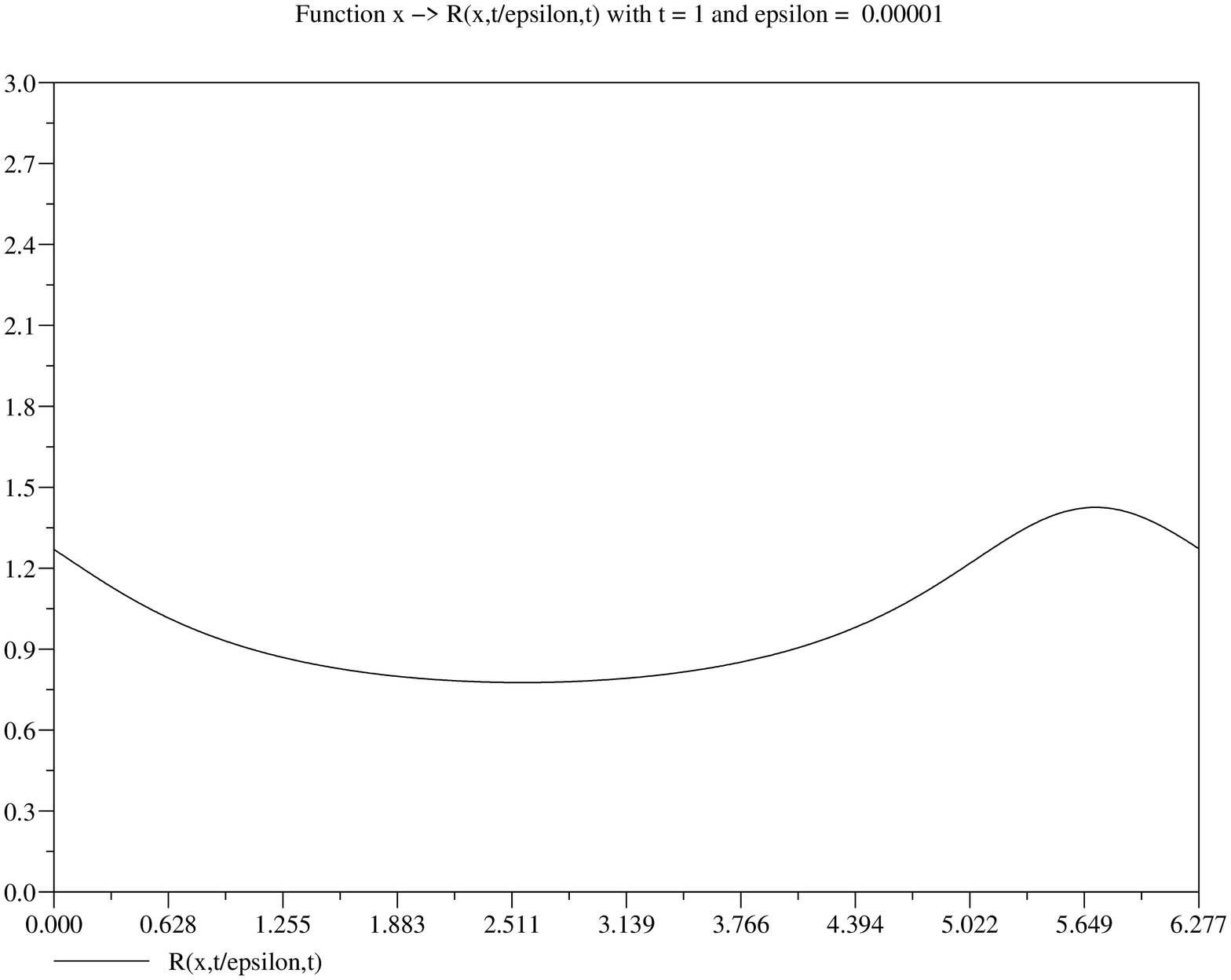} \\
\vspace{-0.5cm}
\includegraphics[scale=0.26,angle=270]{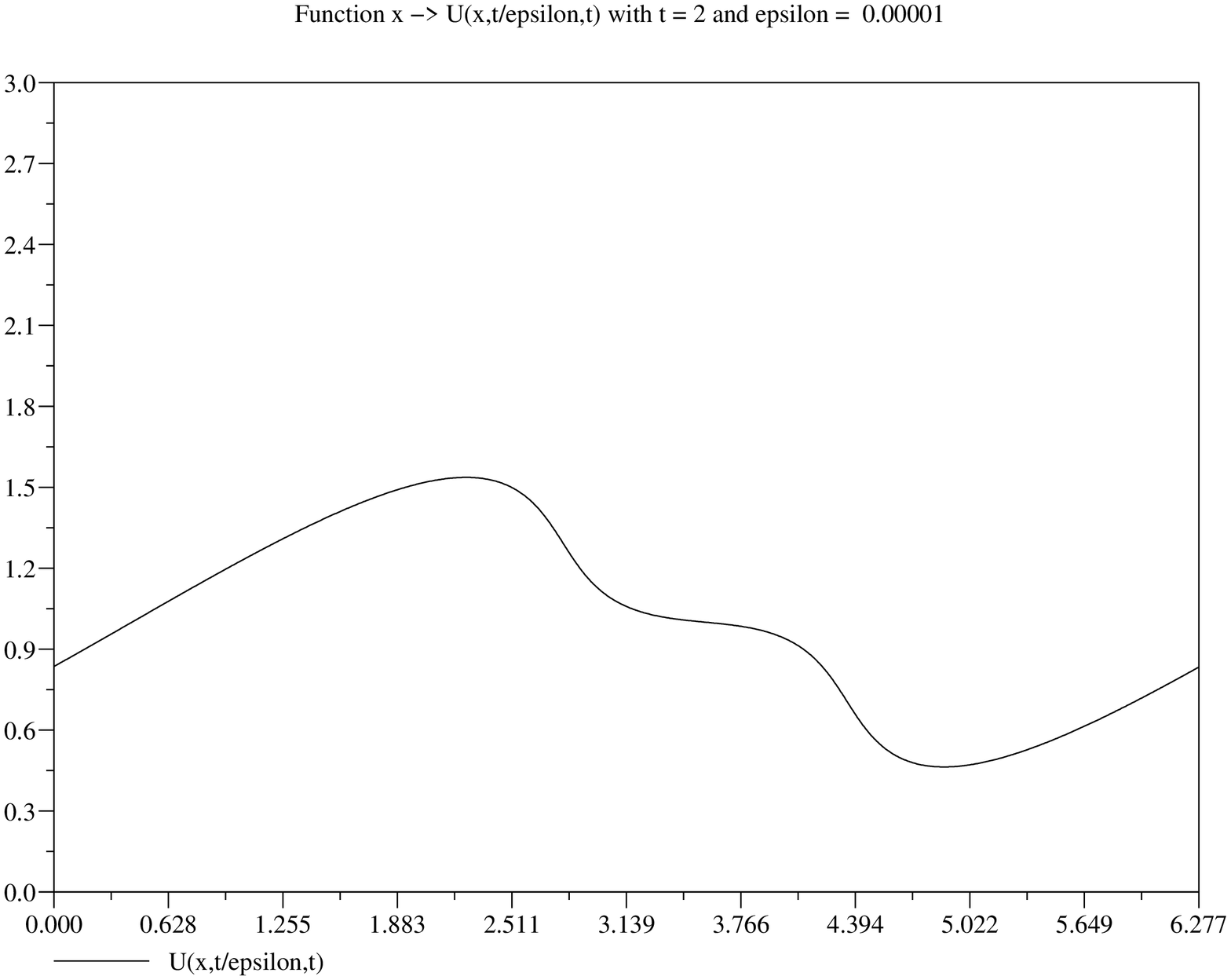} \hspace{-1.15cm}
\includegraphics[scale=0.26,angle=270]{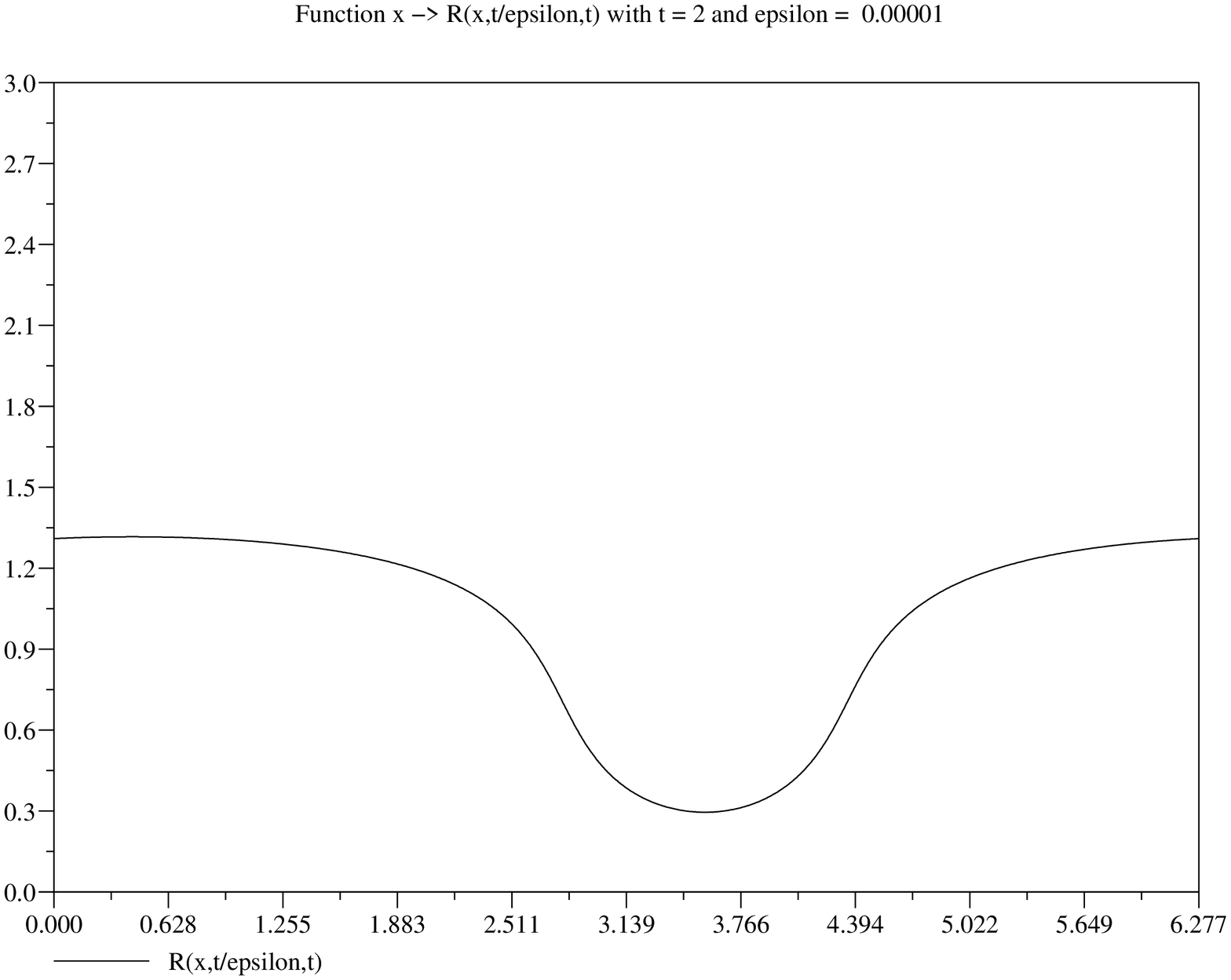} 
\vspace{-0.5cm}
\caption{$U_{h}(x,\frac{t}{\epsilon},t)$ (left) and $R_{h}(x,\frac{t}{\epsilon},t)$ (right) with $\epsilon = 10^{-5}$ at times $t = 0$, $t = 1$, $t = 2$.} 
\end{center}
\begin{center}
\includegraphics[scale=0.26,angle=270]{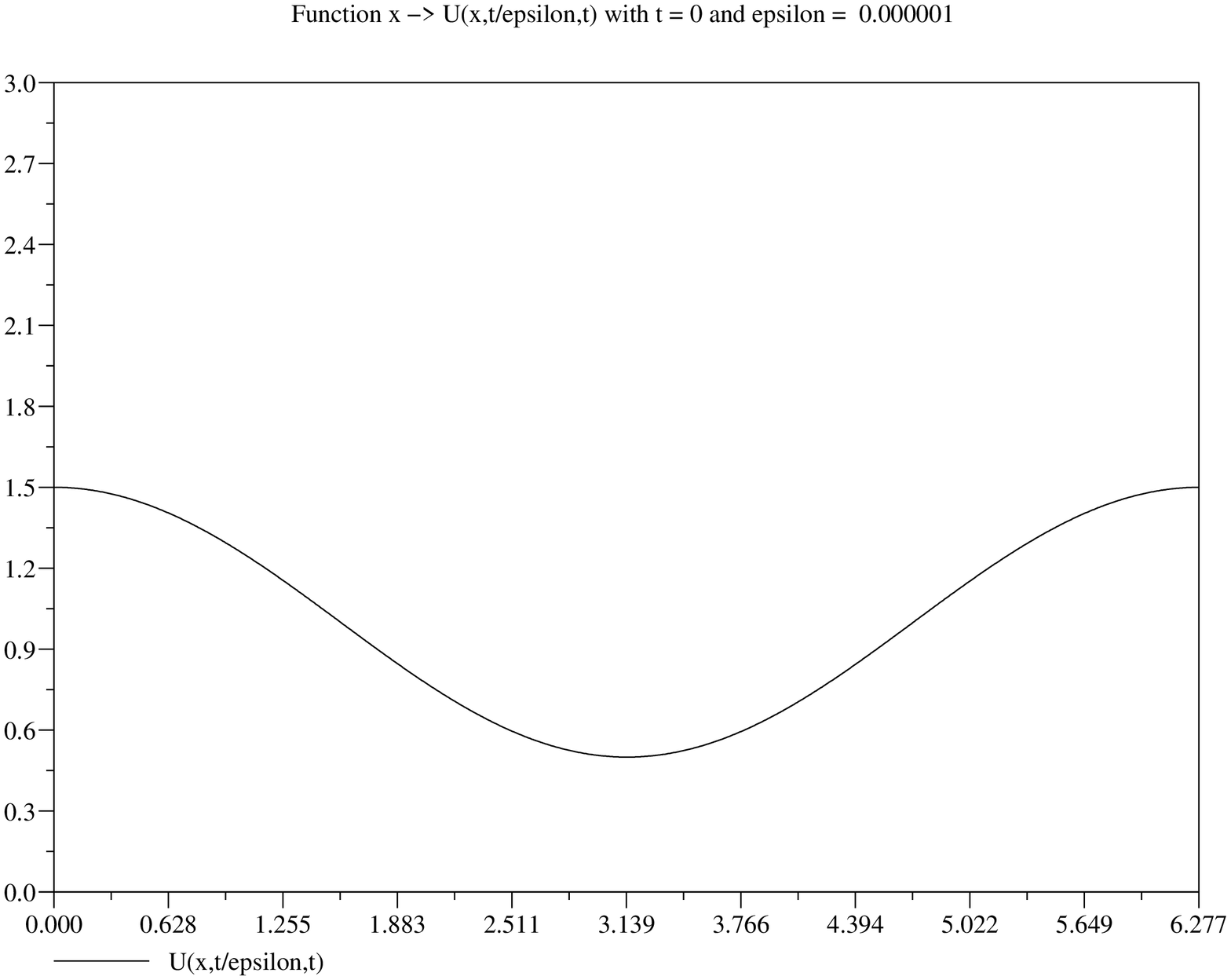} \hspace{-1.15cm}
\includegraphics[scale=0.26,angle=270]{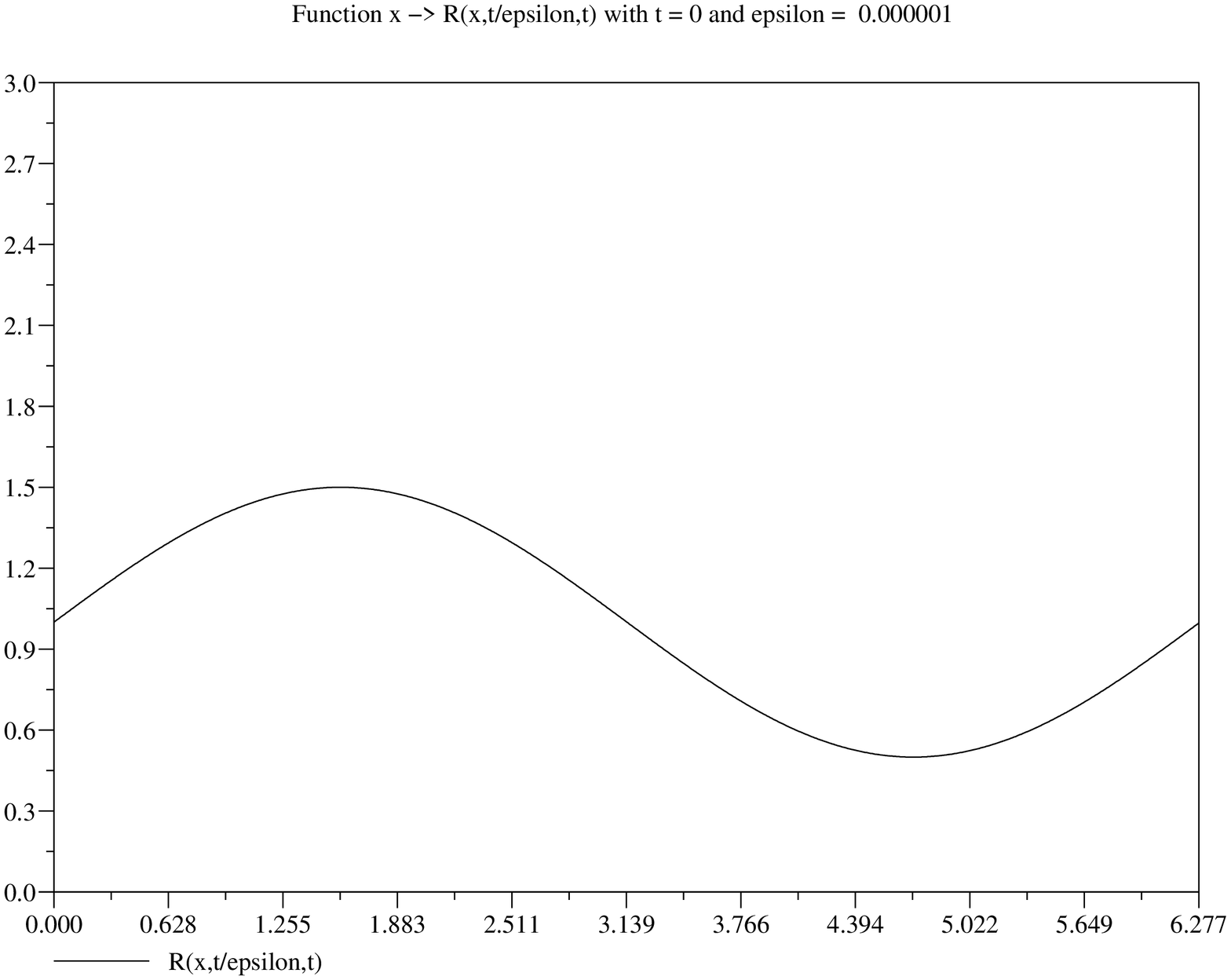} \\
\includegraphics[scale=0.26,angle=270]{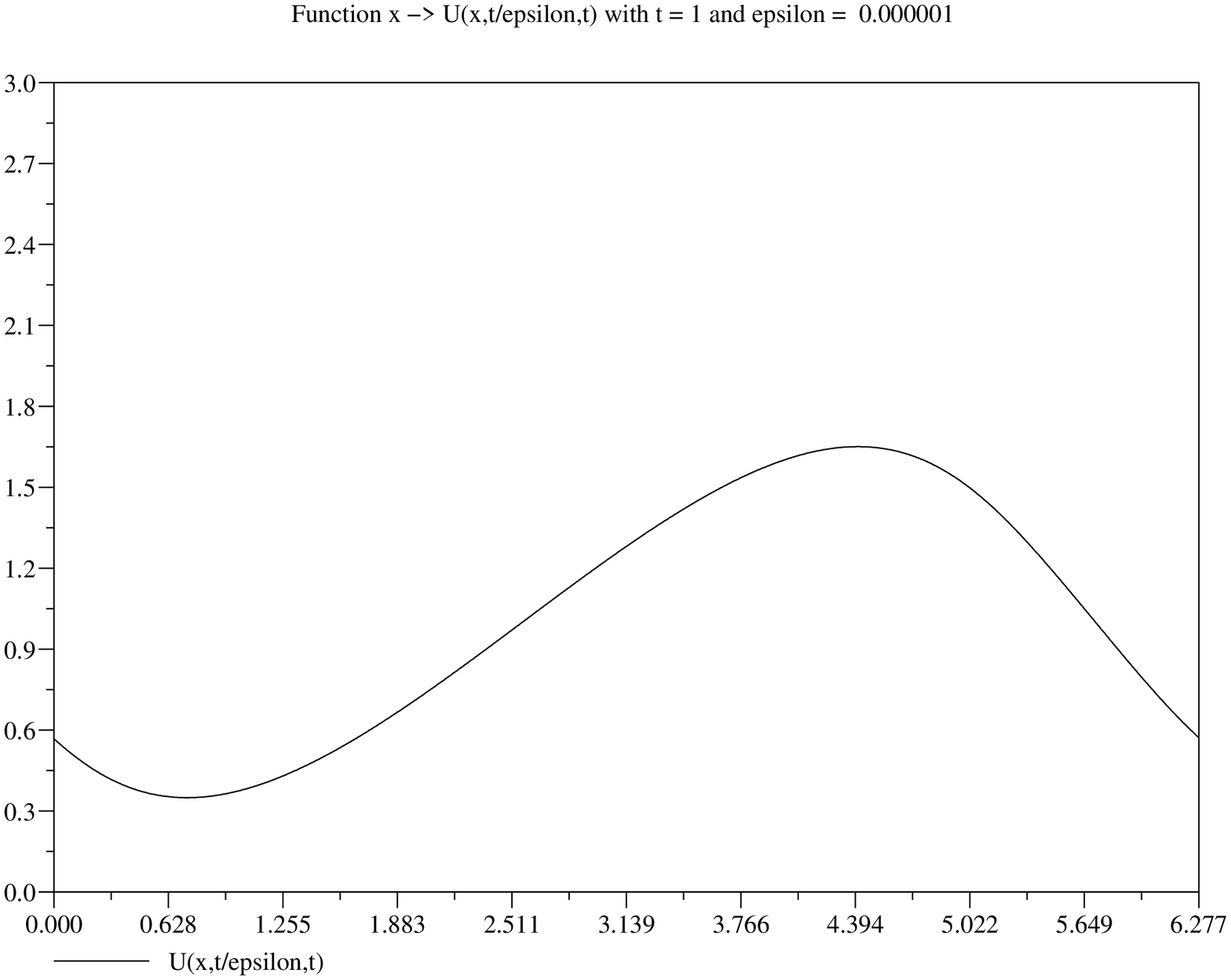} \hspace{-1.15cm}
\includegraphics[scale=0.26,angle=270]{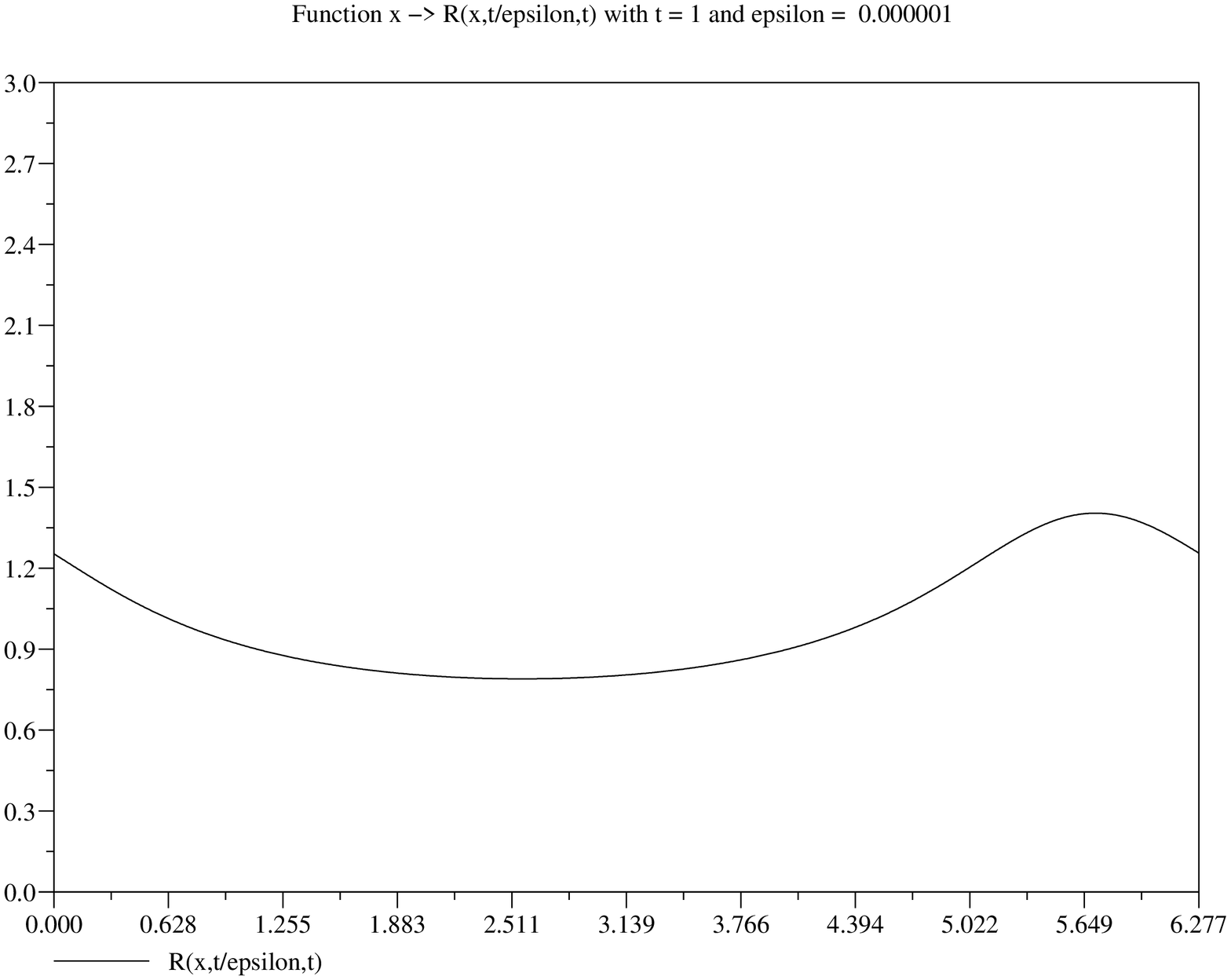}
\includegraphics[scale=0.26,angle=270]{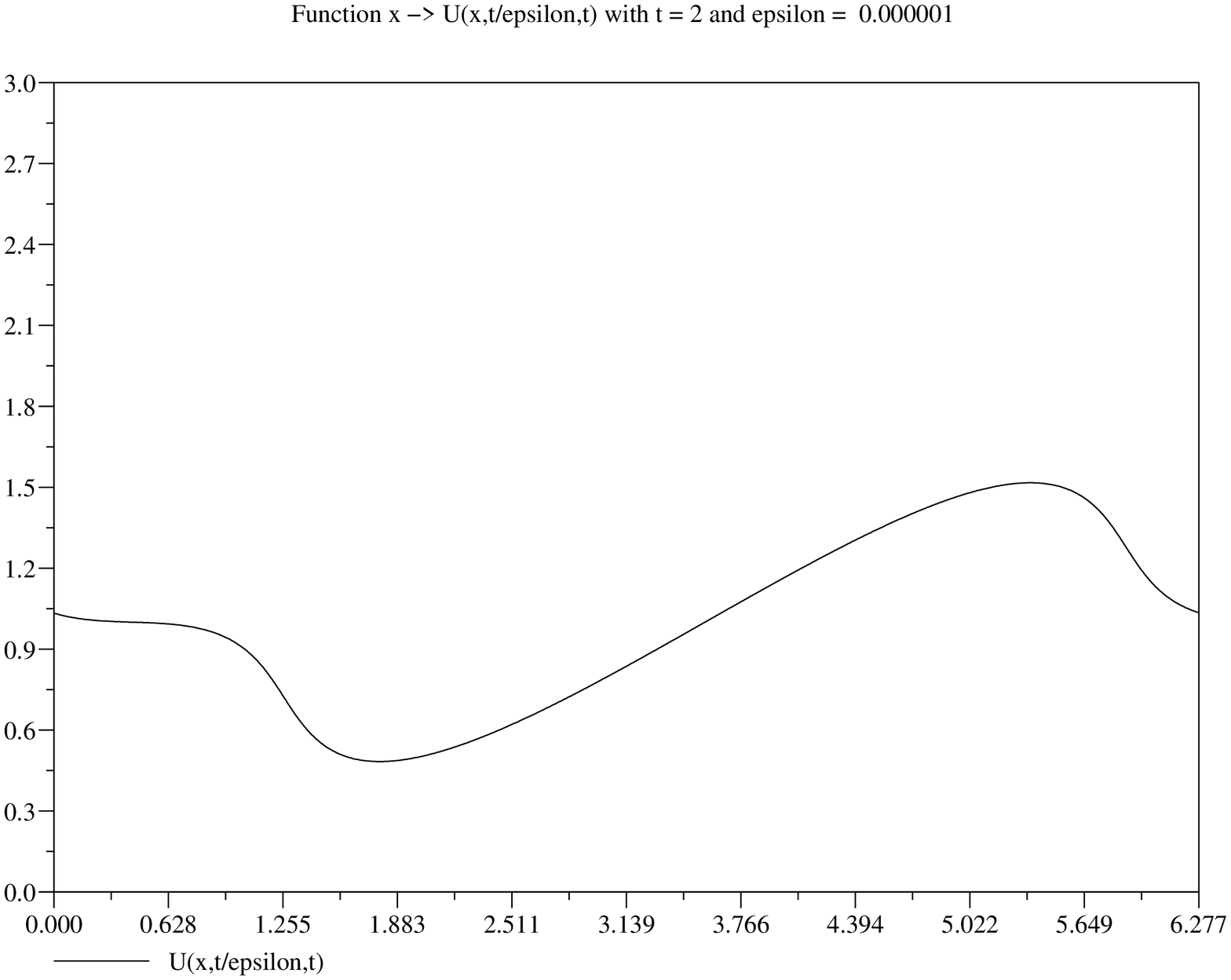} \hspace{-1.15cm}
\includegraphics[scale=0.26,angle=270]{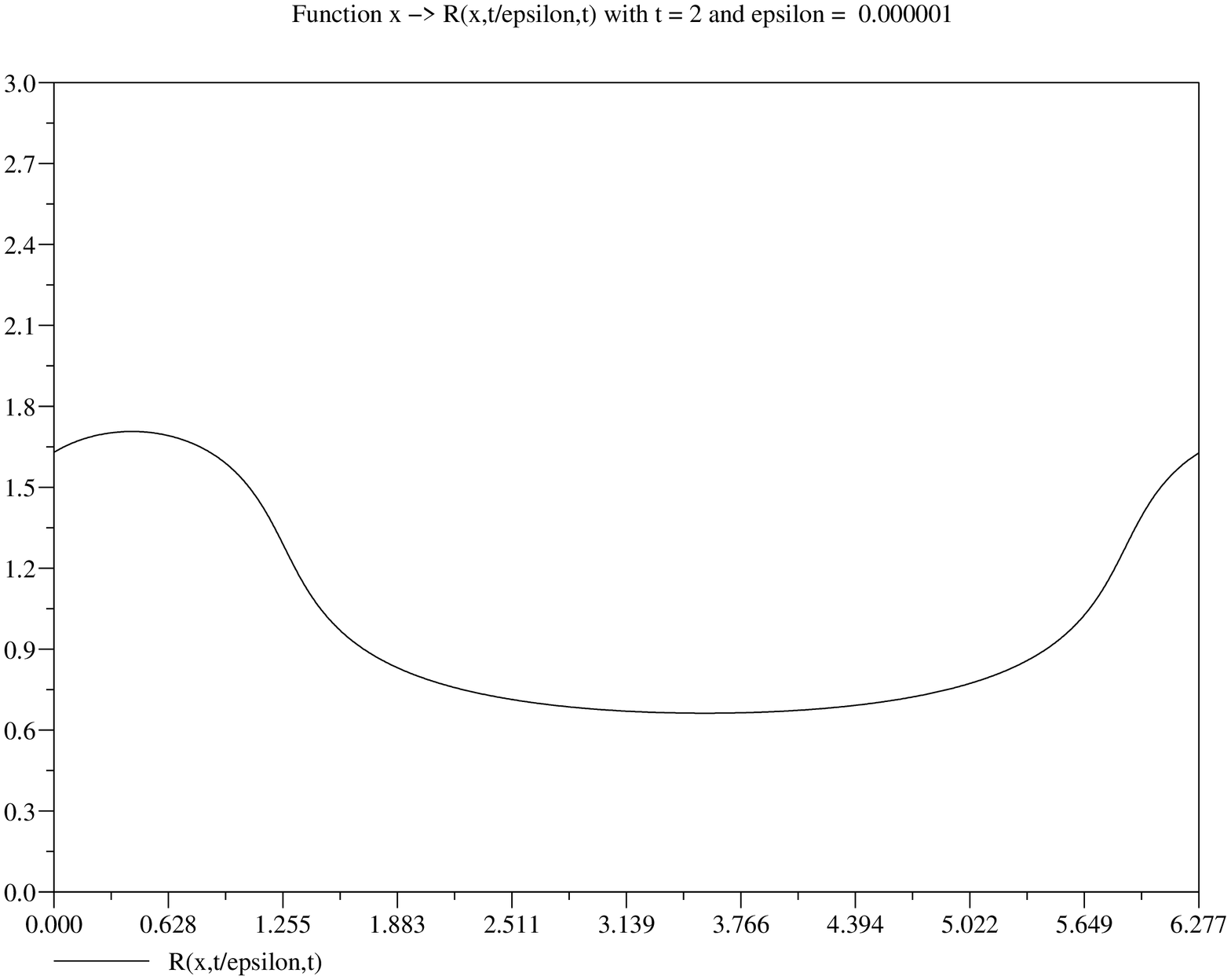}
\caption{$U_{h}(x,\frac{t}{\epsilon},t)$ (left) and $R_{h}(x,\frac{t}{\epsilon},t)$ (right) with $\epsilon = 10^{-6}$ at times $t = 0$, $t = 1$, $t = 2$.}
\end{center}

We can notice that these results have been obtained after a few minutes of computation just as described in the table below:

\begin{center}
\begin{tabular}{|c|c|c|}
\hline
Value of $\epsilon$ & Computation of $(F_{h},B_{h})$ & Computation of $(U_{h},R_{h})$ \\
\hline
$10^{-4}$ & 3 m 45 s 71'' & 42 s 84'' \\
\hline
$10^{-5}$ & 3 m 45 s 95'' & 42 s 87'' \\
\hline
$10^{-6}$ & 3 m 45 s 75'' & 42 s 98'' \\
\hline
\end{tabular}
\begin{table}[ht]
\caption{CPU time costs.}
\end{table}
\end{center}

Furthermore, if we continue the simulation to final time $T \approx 3.2$ for $\epsilon$ ranging in $\{10^{-1},10^{-2},10^{-3},10^{-4},10^{-5},10^{-6}\}$, we remark that discontinuities appear at time $t \approx 3.11$ independently of $\epsilon$ and that, for $\epsilon$ ranging in $\{ 0.01, 0.03, 0.05, 0.07, 0.1\}$, these discontinuities also appear in the solution $(u^{\epsilon},\rho^{\epsilon})$ of (\ref{Grenier}) at the same time.

\section{Conclusion}

\indent In this paper, we have developped a two-scale numerical method for the weakly compressible 1D Euler equations. By using two-scale convergence tools from Nguetseng\cite{General_convergence} and Allaire\cite{Homogenization}, we have developped a new model independent of the Mach number $\epsilon$. Then, using a finite volume scheme on this model, we have obtained numerical approximations that converge to the solutions of the weakly compressible 1D Euler equation when we refine the mesh on $\T^{1} \times [0,T]$ and when $\epsilon \to 0$. \\
\indent Furthermore, we have proved that this two-scale numerical method is first order accurate in $\epsilon$. One the other hand, numerical results have confirmed this accuracy in $\epsilon$ and also have proved that, even with a very small Mach number, the two-scale numerical method allows us to simulate experiments with a reasonable CPU time cost. \\
\indent Motived by the behaviour of the two-scale method, we are now ready to develop such a method on a model which describes a plasma submitted to a strong magnetic field such as the Vlasov-Poisson model established in Fr\'enod and Sonnendr\"ucker \cite{Finite_Larmor_radius}. Furthermore, we only studied the behaviour of the two-scale numerical method in front of smooth initial data: it can be interesting to study the behaviour of the method on non-smooth initial data. Concerning the periodicity, we have built our numerical method in a periodic framework for mathematical reasons. However, it must be possible to extend this method to almost periodic framework in order to study problems like the noise generated by the blades of a turbine or a fan. If we go further, it can be interesting to improve the method in order to use it for studying the noise generated by multidimensional structures like cooling systems for computers or air conditioning, even if the links between time scales and space scales are not clear in this framework.

\end{document}